\documentclass[a4paper,10pt,12pt]{article}
\usepackage{amsmath, amssymb, latexsym, amscd, amsthm,amsfonts,amstext}
\usepackage[mathscr]{eucal}
\usepackage{graphicx}
\usepackage{subfig}

\numberwithin{equation}{section}
\input{tcilatex}
\makeatletter
\renewcommand\tableofcontents{    \@starttoc{toc}}
\makeatother

\begin{document}

\title{Carleman weight functions for a globally convergent numerical method
for ill-posed Cauchy problems for some quasilinear PDEs}
\author{Anatoly B.\ Bakushinskii$^{\ast }$, Michael V. Klibanov$^{\circ
\diamond }$ and Nikolaj A. Koshev$^{\ast \ast }$ \\
%EndAName
\\
$^{\ast }$Institute for Systems Analysis\\
of Russian Academy of Sciences, \\
60 Oktober Anniversary avenue, 9, 117312, Moscow , Russia\\
$^{\circ }$Department of Mathematics \& Statistics, University of North
Carolina \\
at Charlotte, Charlotte, NC 28223, USA\\
$^{\ast \ast }$Institute of Computational Mathematics\\
University of S\~{a}n Paulo, S\~{a}o Carlos, SP 13566-590, Brazil\\
$^{^{\diamond }}$The corresponding author\\
Emails: \texttt{\ bakush@isa.ru, }mklibanv@uncc.edu,
nikolay.koshev@gmail.com }
\date{}
\maketitle

\begin{abstract}
In a series of publications of the second author, including some with
coauthors, globally strictly convex Tikhonov-like functionals were
constructed for some nonlinear ill-posed problems. The main element of such
a functional is the presence of the Carleman Weight Function. Compared with
previous publications, the main novelty of this paper is that the existence
of the regularized solution (i.e. the minimizer) is proved rather than
assumed. The method works for both ill-posed Cauchy problems for some
quasilinear PDEs of the second order and for some Coefficient Inverse
Problems. However, to simplify the presentation, we focus here only on
ill-posed Cauchy problems. Along with the theory, numerical results are
presented for the case of a 1-D quasiliear parabolic PDE with the lateral
Cauchy data given on one edge of the interval (0,1).
\end{abstract}

% -------------------------------------------------------


\section{Introduction}

\label{sec:1}

In this paper we eliminate a restrictive assumption, which was imposed in
the work \cite{Kl1} of the second author. More precisely, the existence of a
minimizer of a weighted Tikhonov functional is proved here rather than
assumed as in \cite{Kl1}. Although similar assumptions of works of the
second author \cite{BKconv,Klib97,Kpar,KNT,KK,KKLY} concerning Coefficient
Inverse Problems (CIPs) can also be eliminated the same way, we are not
doing this here for brevity. In addition to the theory, we present results
of some numerical experiments in which we solve an ill-posed problem for a
1-D quasilinear parabolic equation with the lateral Cauchy data.\ In this
problem, which is also called \emph{side Cauchy problem}, both Dirichlet and
Neumann boundary conditions are given on one edge of the interval $x\in
(0,1) $, the initial condition is unknown and it is required to find the
solution of that equation inside of that interval.

Side Cauchy problems for quasilinear parabolic equations have applications
in processes involving high temperatures \cite{Alif1,Alif2}. In such a
process one can measure both the temperature and the heat flux on one side
of the boundary. However, it is impossible to measure these quantities on
the rest of the boundary. Still, one is required to compute the temperature
in at least a part of the domain of interest. The underlying PDE, which
governs the process of the propagation of this temperature, is a parabolic
PDE. This equation is quasilinear rather then linear because of high
temperatures. The second application is in the glaciology \cite{Koz,Col}. In
this case one is interested in the velocity of a glacier on its bottom side,
which is deeply under the surface of the water. This is the so-called basal
velocity. However, it is impossible to measure that velocity deeply under
the surface of the water. On the other hand, it is possible to measure that
velocity and its normal derivative on the part of the water surface. So, the
idea is to use these measurements to figure out basal velocity. Thus, we
come up with the Cauchy problem for a quasilinear elliptic equation \cite%
{Koz}.

It is well known that the phenomena of multiple local minima and ravines
plagues least squares Tikhonov functionals for nonlinear ill-posed problems,
see, e.g. numerical examples in \cite{KT,Scales}. Therefore, the convergence
of an optimization method for such a functional can be guaranteed only if
its starting point is located in a sufficiently small neighborhood of the
exact solution, i.e. this is \emph{local convergence}. On the other hand, we
call a numerical method for an ill-posed problem \emph{globally convergent}
if there is a theorem, which guarantees that this method delivers at least
one point in a sufficiently small neighborhood of the exact solution without
any advanced knowledge of this neighborhood \cite{BK,Kl1}.

In a series of publications of the second author, including some with
coauthors, starting in 1997 \cite{Klib97,Kpar} and with the recently renewed
interest in \cite{BKconv,KNT,KK} special Tikhonov-like cost functionals $%
J_{\lambda }$ were constructed for CIPs. In particular, some numerical
examples are presented in \cite{KNT}. In \cite{Kl1} this idea was extended
to ill-posed Cauchy problems for quasilinear PDEs of the second order.
Numerical studies of the idea of \cite{Kl1} can be found in \cite{KKLY}. The
key element of each of these functionals is the presence of the Carleman
Weight Function (CWF), i.e. the function which is involved in the Carleman
estimate for the principal part of the corresponding Partial Differential
Operator. The main theorem of each of these works claims that, given a
reasonable bounded set $\Phi $ of an arbitrary diameter in a reasonable
space $H^{k},$ one can choose the parameter $\lambda _{0}=\lambda _{0}\left(
\Phi \right) >1$ of the CWF, depending on $\Phi $, such that for all $%
\lambda \geq \lambda _{0}$ the functional $J_{\lambda }$ is strictly convex
on the set $\Phi $.

The strict convexity implies the absence of multiple local minima and
ravines. Next, the existence of the minimizer of $J_{\lambda }$ on $\Phi $
was assumed. Using this assumption, it was proven that the gradient method
of the minimization of $J_{\lambda }$ converges to that minimizer starting
from an arbitrary point of $\Phi $, provided that all points obtained via
iterations of the gradient method belong to $\Phi .$ Furthermore, it was
established that the distance between that minimizer and the exact solution
of the corresponding inverse problem is small as long as the noise in the
data is small.\ In other words, convergence of regularized solutions was
established. Thus, the above means the \emph{global convergence} of the
gradient method to the exact solution. Still, the assumptions about the
existence of the minimizer on the set $\Phi $ and that all points of the
sequence of the gradient method belong to $\Phi $ are restrictive ones.

In this paper we remove these assumptions via bringing in some ideas of the
convex analysis. To simplify the presentation, we focus here on ill-posed
Cauchy problems for quasilinear PDEs of the second order, i.e. we generalize
results of \cite{Kl1}. We point out, however, that very similar
generalizations can be done for coefficient inverse problems, which were
considered in the above cited works \cite{BKconv,Klib97,Kpar,KNT,KK}.

Those results of the convex analysis require us to change the previous
scheme of the method. More precisely, while the previous scheme of \cite%
{Kl1,KKLY} works with non-zero Cauchy data, we now need to have zero Cauchy
data. We obtain them via \textquotedblleft subtracting" the non-zero Cauchy
data from the sought for solution. In addition, we now need to prove the
Lipschitz continuity of the Frech\'{e}t derivative of our cost functional,
which was not done in those previous works. These factors, in turn mean that
proofs of main theorems here are different from their analogs in \cite%
{Kl1,KKLY}. So, we prove the corresponding theorems below.

The idea of applications of Carleman estimates to CIPs was first published
in 1981 in the work \cite{BukhK}. The method of \cite{BukhK} was originally
designed for proofs of uniqueness theorems for CIPs with single measurement
data, see, e.g. some follow up publications in \cite{Bukh,K92,KT}. There is
now a large number of publications of different authors discussing the idea
of \cite{BukhK}. Since this is not a survey of that method, we cite here
only a few of them \cite{Bukh,Is,K92,Trig}. Surveys of works on the method
of \cite{BukhK} can be found in \cite{Ksurvey,Y}, also, see sections 1.10
and 1.11 of the book \cite{BK}.

In section 2 we present required facts from the convex analysis. In section
3 we present the general scheme of our numerical method for ill-posed Cauchy
problems for quasilinear PDEs of the second order. We also formulate
theorems in section 3. In sections 4-7 we prove those theorems. In section 8
we specify PDEs of the second order for which our technique is applicable.
In section 9 we present numerical results. Summary is presented in section
10.

\section{Some facts of the convex analysis}

\label{sec:2}

Results of this section are known and can be found in chapters 4 and 5 of
the book of Vasiliev \cite{Vas}. Still, we prove below Lemmata 2.1, 2.3 and
Theorem 2.1 for the convenience of the reader. Even though all results of
this section are formulated for a strictly convex functional, some of them
are valid under less restrictive condition, which we do not list here for
brevity.

Let $H$ be a Hilbert space of real valued functions. Below in this section $%
\left\Vert \cdot \right\Vert $ and $\left( ,\right) $ denote the norm and
the scalar product in this space respectively. Let $B\left( R\right)
=\left\{ x\in H:\left\Vert x\right\Vert <R\right\} \subset H$ be the ball of
the radius $R$ with the center at $\left\{ 0\right\} .$ Even though results
of this section can be easily extended to the case when $B\left( R\right) $
is a convex bounded set, we are not doing this here for brevity. Let $\delta
>0$ be a sufficiently small number. Let $J:B\left( R+\delta \right)
\rightarrow \mathbb{R}$ be a functional, which has Frech\'{e}t derivative $%
J^{\prime }\left( x\right) ,\forall x\in \overline{B}\left( R\right) .$
Below we sometimes denote the action of the functional $J^{\prime }\left(
x\right) $ at the point $x$ on any element $h\in H$ as $J^{\prime }\left(
x\right) \left( h\right) .$ But sometimes we also denote this action as $%
\left( J^{\prime }\left( x\right) ,h\right) .$ This difference will not lead
to a misunderstanding. The Frech\'{e}t derivative $J^{\prime }\left(
x\right) $ at a point $x\in \left\{ \left\Vert x\right\Vert =R\right\} $ is
understood as%
\begin{equation*}
J\left( y\right) -J\left( x\right) =J^{\prime }\left( x\right) \left(
y-x\right) +o\left( \left\Vert x-y\right\Vert \right) ,\left\Vert
x-y\right\Vert \rightarrow 0,y\in B\left( R+\delta \right) .
\end{equation*}%
We assume that this Frech\'{e}t derivative satisfies the Lipschitz
continuity condition,%
\begin{equation}
\left\Vert J^{\prime }\left( x\right) -J^{\prime }\left( y\right)
\right\Vert \leq L\left\Vert x-y\right\Vert ,\forall x,y\in \overline{B}%
\left( R\right) ,  \label{8.1}
\end{equation}%
with a certain constant $L>0$. In addition, we assume that the functional $%
J\left( x\right) $ is strictly convex on the set $\overline{B}\left(
R\right) ,$%
\begin{equation}
J\left( y\right) -J\left( x\right) -J^{\prime }\left( x\right) \left(
y-x\right) \geq \varkappa \left\Vert x-y\right\Vert ^{2},\forall x,y\in 
\overline{B}\left( R\right) ,  \label{8.2}
\end{equation}%
where $\varkappa =const.>0.$ The strict convexity of $J\left( x\right) $ on $%
\overline{B}\left( R\right) $ implies 
\begin{equation}
\left( J^{\prime }\left( x\right) -J^{\prime }\left( y\right) ,x-y\right)
\geq 2\varkappa \left\Vert x-y\right\Vert ^{2},\forall x,y\in \overline{B}%
\left( R\right) .  \label{8.20}
\end{equation}

\textbf{Lemma 2.1}. \emph{A point }$x_{\min }\in \overline{B}\left( R\right) 
$\emph{\ is a point of a relative minimum of the functional }$J\left(
x\right) $\emph{\ on the set }$\overline{B}\left( R\right) $ \emph{if and
only if} \emph{\ \ }%
\begin{equation}
\left( J^{\prime }\left( x_{\min }\right) ,x_{\min }-y\right) \leq 0,\forall
y\in \overline{B}\left( R\right) .  \label{8.3}
\end{equation}%
\emph{If a point }$x_{\min }\in \overline{B}\left( R\right) $\emph{\ is a
point of a relative minimum of the functional }$J\left( x\right) $\emph{\ on
the set }$\overline{B}\left( R\right) ,$\emph{\ then this point is unique
and it is, therefore, the point of the unique global minimum of }$J\left(
x\right) $ \emph{on the set }$\overline{B}\left( R\right) .$

Note that if $x_{\min }$ is an interior point of $B\left( R\right) $, then
in (\ref{8.3}) \textquotedblleft $\leq $" must be replaced with
\textquotedblleft $=$" and the assertion of this lemma becomes obvious.
However this assertion is not immediately obvious if $x_{\min }$ belongs to
the boundary of the closed ball $\overline{B}\left( R\right) .$

\textbf{Proof}. Suppose that $x_{\min }$ is a point of a relative minimum of 
$J\left( x\right) $ on $\overline{B}\left( R\right) .$ Assume to the
contrary: that there exists a point $y\in \overline{B}\left( R\right) $ such
that $\left( J^{\prime }\left( x_{\min }\right) ,x_{\min }-y\right) >0.$ Let 
$h=y-x_{\min }.$ Then 
\begin{equation}
\left( J^{\prime }\left( x_{\min }\right) ,\xi h\right) <0,\forall \xi >0
\label{8.4}
\end{equation}%
for any number $\xi >0.$\ Since the set $\overline{B}\left( R\right) $ is
convex, then $\left\{ x_{\min }+\xi h,\xi \in \left[ 0,1\right] \right\}
\subset \overline{B}\left( R\right) .$ We have%
\begin{equation}
J\left( x_{\min }+\xi h\right) =J\left( x_{\min }\right) +\xi \left[ \left(
J^{\prime }\left( x_{\min }\right) ,h\right) +o\left( 1\right) \right] ,\xi
\rightarrow 0^{+}.  \label{8.40}
\end{equation}%
By (\ref{8.4}) $\left( J^{\prime }\left( x_{\min }\right) ,h\right) +o\left(
1\right) <0$ for sufficiently small values of $\xi >0.$ Hence, (\ref{8.40})
implies that $J\left( x_{\min }+\xi h\right) <J\left( x_{\min }\right) $ for
sufficiently small $\xi .$ The latter contradicts the assumption that $%
x_{\min }$ is a point of a relative minimum of the functional $J\left(
x\right) $ on the set $\overline{B}\left( R\right) .$

Assume now the reverse: that the inequality (\ref{8.3}) is valid for a
certain point $x_{\min }\in \overline{B}\left( R\right) .$ We prove below
that $x_{\min }$ is a point of \ a relative minimum of the functional $%
J\left( x\right) $ on the set $\overline{B}\left( R\right) .$ Indeed, let $%
y\in \overline{B}\left( R\right) $ be an arbitrary point and let $y\neq x$.
By (\ref{8.3}) $J^{\prime }\left( x_{\min }\right) \left( y-x_{\min }\right)
\geq 0.$ Hence, (\ref{8.2}) implies that%
\begin{equation}
J\left( y\right) \geq J\left( x_{\min }\right) +J^{\prime }\left( x_{\min
}\right) \left( y-x_{\min }\right) +\varkappa \left\Vert x-y\right\Vert
^{2}>J\left( x_{\min }\right) .  \label{8.5}
\end{equation}%
Hence, the functional $J\left( x\right) $ attains its minimal value at $%
x=x_{\min }$.\ Hence, $x_{\min }$ is indeed the point of a relative minimum
of the functional $J\left( x\right) $ on the set $\overline{B}\left(
R\right) $.

We now prove uniqueness of the point of a relative minimum. Indeed, assume
that there are two points $x_{\min }$ and $y_{\min }$ of relative minima of
the functional $J\left( x\right) $ on the set $\overline{B}\left( R\right) $%
. We have%
\begin{equation}
J\left( y_{\min }\right) -J\left( x_{\min }\right) -J^{\prime }\left(
x_{\min }\right) \left( y_{\min }-x_{\min }\right) \geq \varkappa \left\Vert
x_{\min }-y_{\min }\right\Vert ^{2},  \label{8.51}
\end{equation}%
\begin{equation}
J\left( x_{\min }\right) -J\left( y_{\min }\right) -J^{\prime }\left(
y_{\min }\right) \left( x_{\min }-y_{\min }\right) \geq \varkappa \left\Vert
x_{\min }-y_{\min }\right\Vert ^{2}.  \label{8.52}
\end{equation}%
Summing up (\ref{8.51}) and (\ref{8.52}), we obtain%
\begin{equation}
-J^{\prime }\left( x_{\min }\right) \left( y_{\min }-x_{\min }\right)
-J^{\prime }\left( y_{\min }\right) \left( x_{\min }-y_{\min }\right) \geq
2\varkappa \left\Vert x_{\min }-y_{\min }\right\Vert ^{2}.  \label{8.53}
\end{equation}%
However, by (\ref{8.3}) 
\begin{equation}
-J^{\prime }\left( x_{\min }\right) \left( y_{\min }-x_{\min }\right)
-J^{\prime }\left( y_{\min }\right) \left( x_{\min }-y_{\min }\right) \leq 0.
\label{8.54}
\end{equation}%
Hence, (\ref{8.53}) and (\ref{8.54}) imply that $x_{\min }=y_{\min }.$ $%
\square $

Let $y\in H$ be an arbitrary point. The point $\overline{y}$ is called
projection of the point $y$ on the set $\overline{B}\left( R\right) $ if%
\begin{equation*}
\left\Vert y-\overline{y}\right\Vert =\inf_{v\in \overline{B}\left( R\right)
}\left\Vert y-v\right\Vert .
\end{equation*}

\textbf{Lemma 2.2}. \emph{Each point }$y\in H$\emph{\ has unique projection }%
$\overline{y}$\emph{\ on the set }$\overline{B}\left( R\right) .$\emph{\
Furthermore, the point }$\overline{y}\in \overline{B}\left( R\right) $\emph{%
\ is the projection of the point }$y$\emph{\ on the set }$\overline{B}\left(
R\right) $\emph{\ if and only if }%
\begin{equation}
\left( \overline{y}-y,v-\overline{y}\right) \geq 0,\forall v\in \overline{B}%
\left( R\right) .  \label{8.7}
\end{equation}

For the proof of this lemma we refer to theorem 1 of \S 4 of chapter 4 of 
\cite{Vas}. Denote the projection operator of the space $H$ on the set $%
\overline{B}\left( R\right) $ as $P_{\overline{B}\left( R\right)
}:H\rightarrow \overline{B}\left( R\right) .$ Then (see theorem 2 of \S 4 of
chapter 4 of \cite{Vas}) 
\begin{equation}
\left\Vert P_{\overline{B}\left( R\right) }\left( u\right) -P_{\overline{B}%
\left( R\right) }\left( v\right) \right\Vert \leq \left\Vert u-v\right\Vert
,\forall u,v\in H.  \label{8.8}
\end{equation}

\textbf{Lemma 2.3}. \emph{The point }$x_{\min }\in \overline{B}\left(
R\right) $\emph{\ is the point of the unique global minimum of the
functional }$J\left( x\right) $\emph{\ on the set }$\overline{B}\left(
R\right) $\emph{\ if and only if there exits a number }$\gamma >0$ \emph{%
such that} 
\begin{equation}
x_{\min }=P_{\overline{B}\left( R\right) }\left( x_{\min }-\gamma J^{\prime
}\left( x_{\min }\right) \right) .  \label{8.9}
\end{equation}%
\emph{If (\ref{8.9}) is valid for one number }$\gamma ,$\emph{\ then it is
also valid for all }$\gamma >0.$

\textbf{Proof}. Uniqueness of the global minimum, if it exists, and the
absence of other relative minima, follows from Lemma 2.1. By (\ref{8.7})
equality (\ref{8.9}) is equivalent with%
\begin{equation}
\left( x_{\min }-\left( x_{\min }-\gamma J^{\prime }\left( x_{\min }\right)
\right) ,v-x_{\min }\right) \geq 0,\forall v\in \overline{B}\left( R\right) .
\label{8.10}
\end{equation}%
Since $\gamma >0,$ then (\ref{8.10}) implies that $\left( J^{\prime }\left(
x_{\min }\right) ,x_{\min }-v\right) \leq 0,\forall v\in \overline{B}\left(
R\right) ,$ which is exactly (\ref{8.3}). The rest follows immediately from
Lemma 2.1. $\square $

Consider now the gradient projection method to find the minimum of the
functional $J\left( x\right) $ on the set $\overline{B}\left( R\right) .$
Let $x_{0}\in \overline{B}\left( R\right) $ be an arbitrary point. We
construct the following sequence%
\begin{equation}
x_{n+1}=P_{\overline{B}\left( R\right) }\left( x_{n}-\gamma J^{\prime
}\left( x_{n}\right) \right) ,n=0,1,2,...  \label{8.11}
\end{equation}

\textbf{Theorem 2.1}. \emph{Assume that the functional }$J\left( x\right) $%
\emph{\ is strictly convex on the closed ball }$\overline{B}\left( R\right) $%
\emph{\ and let condition (\ref{8.1}) holds.\ Then there exists unique point
of the relative minimum }$x_{\min }$\emph{\ of this functional on the set }$%
\overline{B}\left( R\right) .$\emph{\ Furthermore, }$x_{\min }$ \emph{is the
unique point of the global minimum of }$J\left( x\right) $\emph{\ on }$%
\overline{B}\left( R\right) .$\emph{\ Let }$L$\emph{\ and }$\gamma $\emph{\
be numbers in (\ref{8.1}) and (\ref{8.2}) respectively and let }$\gamma \in
\left( 0,L\right] .$\emph{\ Let the number }$\gamma $\emph{\ in (\ref{8.11})
be so small that }%
\begin{equation}
0<\gamma <2\varkappa L^{-2}.  \label{8.12}
\end{equation}%
\emph{Let} $q\left( \alpha \right) =\left( 1-2\gamma \varkappa +\alpha
^{2}L^{2}\right) ^{1/2}.$ \emph{Then the sequence (\ref{8.11}) converges to
the point }$x_{\min }$\emph{\ and }%
\begin{equation}
\left\Vert x_{n}-x_{\min }\right\Vert \leq q^{n}\left( \gamma \right)
\left\Vert x_{0}-x_{\min }\right\Vert .  \label{8.13}
\end{equation}

\textbf{Proof}. We note first that by (\ref{8.12}) the number $q\left(
\gamma \right) \in \left( 0,1\right) .$ The idea of the proof is to show
that the operator in the right hand side of (\ref{8.11}) is contraction
mapping, as long as (\ref{8.12})\emph{\ }holds. Denote $D\left( x\right) =P_{%
\overline{B}\left( R\right) }\left( x-\gamma J^{\prime }\left( x\right)
\right) ,x\in \overline{B}\left( R\right) .$ Then the operator$\ D:\overline{%
B}\left( R\right) \rightarrow \overline{B}\left( R\right) .$ Let $x$ and $y$
be two arbitrary points of $\overline{B}\left( R\right) .$ Using (\ref{8.8}%
), we obtain%
\begin{equation}
\begin{array}{c}
\left\Vert D\left( x\right) -D\left( y\right) \right\Vert ^{2}\leq
\left\Vert \left( x-\gamma J^{\prime }\left( x\right) \right) -\left(
y-\gamma J^{\prime }\left( y\right) \right) \right\Vert ^{2} \\ 
=\left\Vert \left( x-y\right) -\gamma \left( J^{\prime }\left( x\right)
-J^{\prime }\left( y\right) \right) \right\Vert ^{2} \\ 
=\left\Vert x-y\right\Vert ^{2}+\gamma ^{2}\left\Vert J^{\prime }\left(
x\right) -J^{\prime }\left( y\right) \right\Vert ^{2}-2\gamma \left(
J^{\prime }\left( x\right) -J^{\prime }\left( y\right) ,x-y\right) .%
\end{array}
\label{8.14}
\end{equation}%
By (\ref{8.1}) $\gamma ^{2}\left\Vert J^{\prime }\left( x\right) -J^{\prime
}\left( y\right) \right\Vert ^{2}\leq \gamma ^{2}L^{2}\left\Vert
x-y\right\Vert ^{2}.$ Next, by (\ref{8.20}) 
\begin{equation*}
-2\gamma \left( J^{\prime }\left( x\right) -J^{\prime }\left( y\right)
,x-y\right) \leq -2\gamma \varkappa \left\Vert x-y\right\Vert ^{2}.
\end{equation*}%
Hence, (\ref{8.14}) leads to%
\begin{equation*}
\left\Vert D\left( x\right) -D\left( y\right) \right\Vert ^{2}\leq \left(
1-2\gamma \varkappa +\gamma ^{2}L^{2}\right) \left\Vert x-y\right\Vert
^{2}=q^{2}\left( \gamma \right) \left\Vert x-y\right\Vert ^{2}.
\end{equation*}%
Hence, the operator $D$ is \ a contraction mapping of the set $\overline{B}%
\left( R\right) .$ The rest of the proof follows immediately from Lemmata
2.1 and 2.3. $\square $

\section{The general scheme of the method}

\label{sec:3}

\subsection{The Cauchy problem}

\label{sec:3.1}

Let $\Omega \subset \mathbb{R}^{n}$ be a bounded domain.\ Let $A$ be a
quasilinear Partial Differential Operator of the second order in $\Omega $
with its linear principal part $A_{0},$ 
\begin{equation}
A\left( u\right) =\sum\limits_{\left\vert \alpha \right\vert =2}a_{\alpha
}\left( x\right) D^{\alpha }u+A_{1}\left( x,\nabla u,u\right) ,
\label{2.100}
\end{equation}%
\begin{equation}
A_{0}u=\sum\limits_{\left\vert \alpha \right\vert =2}a_{\alpha }\left(
x\right) D^{\alpha }u,  \label{2.1001}
\end{equation}%
\begin{equation}
a_{\alpha }\in C^{1}\left( \overline{\Omega }\right) ,  \label{300}
\end{equation}%
\begin{equation}
\text{ }A_{1}\left( x,y\right) \in C^{3}\left\{ \left( x,y\right) :x\in 
\overline{\Omega },y\in \mathbb{R}^{n+1}\right\} .  \label{2.1002}
\end{equation}%
Denote $k=\left[ n/2\right] +2,$ where $\left[ n/2\right] $ is the largest
integer which does not exceed the number $n/2.$ By the embedding theorem 
\begin{equation}
H^{k}\left( \Omega \right) \subset C^{1}\left( \overline{\Omega }\right) 
\text{ and }\left\Vert f\right\Vert _{C^{1}\left( \overline{\Omega }\right)
}\leq C\left\Vert f\right\Vert _{H^{k}\left( \Omega \right) },\forall f\in
H^{k}\left( \Omega \right) ,  \label{2.10005}
\end{equation}%
where the constant $C=C\left( \Omega \right) >0$ depends only on listed
parameters. Let $\Gamma \subseteq \partial \Omega ,\Gamma \in C^{\infty }$
be a part of the boundary of the domain $\Omega .$ We assume that $\Gamma $
is not a part of the characteristic hypersurface of the operator $A_{0}.$

\textbf{Cauchy Problem 1}. \emph{Consider the following Cauchy problem for
the operator }$A,$\emph{\ }%
\begin{equation}
A\left( u\right) =0\text{ in }\Omega ,  \label{2.1003}
\end{equation}%
\begin{equation}
u\mid _{\Gamma }=g_{0}\left( x\right) ,\partial _{n}u\mid _{\Gamma
}=g_{1}\left( x\right) .  \label{2.1004}
\end{equation}%
\emph{Find the solution }$u\in H^{k}\left( \Omega \right) $\emph{\ of the
problem (\ref{2.1003}), (\ref{2.1004}) either in the entire domain }$\Omega $%
\emph{\ or at least in its subdomain.}

The Cauchy-Kowalewski uniqueness theorem is inapplicable here since we do
not impose the analyticity assumption on coefficients $a_{\alpha }\left(
x\right) $ of the principal part $A_{0}$ of the operator $A$ and also since $%
A$ is not a linear operator. Still, Theorem 3.1 guarantees uniqueness of
this problem in the domain $\Omega _{c}$ defined in subsection 3.1.

Suppose that there exists a function $F\in H^{k+1}\left( \Omega \right) $
such that%
\begin{equation}
F\mid _{\Gamma }=g_{0}\left( x\right) ,\partial _{n}F\mid _{\Gamma
}=g_{1}\left( x\right) .  \label{2.1006}
\end{equation}%
Consider the function $v\left( x\right) =u\left( x\right) -F\left( x\right)
. $\ Here is an example of the function $F\left( x\right) $. Suppose that $%
\Omega =\left\{ \left\vert x\right\vert <1\right\} \subset \mathbb{R}^{3}.$
Let $\Gamma =\left\{ \left\vert x\right\vert =1\right\} .$ Assume that
functions $g_{0},g_{1}\in C^{k+1}\left( \Gamma \right) .$ Let the function $%
\chi \left( x\right) \in C^{k+1}\left( \overline{\Omega }\right) $ be such
that%
\begin{equation*}
\chi \left( x\right) =\left\{ 
\begin{array}{c}
1,\left\vert x\right\vert \in \left[ 3/4,1\right] , \\ 
\text{between 0 and 1 for }x\in \left( 1/2,3/4\right) , \\ 
0\text{ for }x\in \left( 0,1/2\right) .%
\end{array}%
\right.
\end{equation*}%
The existence of such functions $\chi \left( x\right) $ is well known from
the Real Analysis course. Then the function $F\left( x\right) $ can be
constructed as $F\left( x\right) =\chi \left( x\right) \left[ g_{0}\left(
x\right) +\left( \left\vert x\right\vert -1\right) g_{1}\left( x\right) %
\right] .$

Define the subspace $H_{0}^{k}\left( \Omega \right) $ of the Hilbert space
of real valued functions $H^{k}\left( \Omega \right) $ as%
\begin{equation*}
H_{0}^{k}\left( \Omega \right) =\left\{ f\in H^{k}\left( \Omega \right)
:f\mid _{\Gamma }=0,\partial _{n}f\mid _{\Gamma }=0\right\} .
\end{equation*}%
Hence, we come up with the following Cauchy problem:

\textbf{Cauchy Problem 2}. \emph{Determine the function }$v\in
H_{0}^{k}\left( \Omega \right) $\emph{\ such that }%
\begin{equation}
A\left( v+F\right) =0\text{ in }\Omega .  \label{2.1007}
\end{equation}

Note that the function $A\left( F\right) \in H^{k-1}\left( \Omega \right) .$
By the embedding theorem, the latter means that $A\left( F\right) \in
C\left( \overline{\Omega }\right) .$ In the realistic case, the Cauchy data $%
g_{0}\left( x\right) ,g_{1}\left( x\right) $ are given with a random noise.
On the other hand, by (\ref{2.1006}) one should have at least the following
smoothness $g_{0}\in H^{k}\left( \Gamma \right) ,g_{1}\in H^{k-1}\left(
\Gamma \right) .$ Hence, a data smoothing procedure might be applied to
these functions in a data pre-processing procedure. A specific form of a
smoothing procedure depends on a specific problem under the consideration.
As a result, one would obtain the Cauchy data with a smooth error. A
smoothing procedure is outside of the scope of this publication. Still, we
work with noisy data in our computations, see section 9.

\subsection{The pointwise Carleman estimate}

\label{sec:3.2}

Let the function $\psi \in C^{\infty }\left( \overline{\Omega }\right) $ and 
$\left\vert \nabla \psi \right\vert \neq 0$ in $\overline{\Omega }.$ For a
number $\alpha >0$ denote 
\begin{equation}
\psi _{\alpha }=\left\{ x\in \overline{\Omega }:\psi \left( x\right) =\alpha
\right\} ,\Omega _{\alpha }=\left\{ x\in \Omega :\psi \left( x\right)
>\alpha \right\} .  \label{2.0}
\end{equation}%
Hence, a part of the boundary $\partial \Omega _{\alpha }$ of the domain $%
\Omega _{\alpha }$ is the level hypersurface $\psi _{\alpha }$ of the
function $\psi .$ We assume that $\Omega _{\alpha }\neq \varnothing .$
Obviously $\Omega _{\omega }\subset \Omega _{\alpha }$ if $\omega >\alpha .$
Choose a sufficiently small number $\varepsilon >0$ such that $\Omega
_{\alpha +2\varepsilon }\neq \varnothing .$ Denote $\Gamma _{\alpha }=\Gamma
\cap \overline{\Omega }_{\alpha }$ and assume that $\Gamma _{\alpha }\neq
\varnothing .$ Hence, the boundary $\partial \Omega _{\alpha }$ of the
domain $\Omega _{\alpha }$ is:%
\begin{equation}
\partial \Omega _{\alpha }=\partial _{1}\Omega _{\alpha }\cup \partial
_{2}\Omega _{\alpha },\text{ }  \label{2.1}
\end{equation}%
\begin{equation}
\partial _{1}\Omega _{\alpha }=\psi _{\alpha },\partial _{2}\Omega _{\alpha
}=\Gamma _{\alpha }.  \label{2.01}
\end{equation}

Let $\lambda >1$ be a large parameter. Consider the function $\varphi
_{\lambda }\left( x\right) ,$%
\begin{equation}
\varphi _{\lambda }\left( x\right) =\exp \left[ \lambda \psi \left( x\right) %
\right] .  \label{2.2}
\end{equation}%
By (\ref{2.1})-(\ref{2.2}) 
\begin{equation}
\min_{\overline{\Omega }_{\alpha }}\varphi _{\lambda }\left( x\right)
=\varphi _{\lambda }\left( x\right) \mid _{\psi _{\alpha }}=e^{\lambda
\alpha }.  \label{2.3}
\end{equation}%
Let%
\begin{equation}
m=\max_{\overline{\Omega }_{\alpha }}\psi \left( x\right) .  \label{2.30}
\end{equation}%
Then%
\begin{equation}
\max_{\overline{\Omega }_{\alpha }}\varphi _{\lambda }\left( x\right)
=e^{\lambda m}.  \label{2.4}
\end{equation}

Assume that the following pointwise estimate is valid for the principal part 
$A_{0}$ of the operator $A:$%
\begin{equation}
\left( A_{0}u\right) ^{2}\varphi _{\lambda }^{2}\left( x\right) \geq
C_{1}\lambda \left( \nabla u\right) ^{2}\varphi _{\lambda }^{2}\left(
x\right) +C_{1}\lambda ^{3}u^{2}\varphi _{\lambda }^{2}\left( x\right) +%
\func{div}U,  \label{2.6}
\end{equation}%
\begin{equation}
U=\left( U_{1},...,U_{n}\right) ,\left\vert U\left( x\right) \right\vert
\leq C_{1}\lambda ^{3}\left[ \left( \nabla u\right) ^{2}+u^{2}\right]
\varphi _{\lambda }^{2}\left( x\right) ,  \label{2.7}
\end{equation}%
\begin{equation}
\forall \lambda \geq \lambda _{0},\forall x\in \Omega _{\alpha },\forall
u\in C^{2}\left( \overline{\Omega }_{\alpha }\right) ,  \label{2.8}
\end{equation}%
where constants $\lambda _{0}=\lambda _{0}\left( A_{0},\Omega \right)
>1,C_{1}=C_{1}\left( A_{0},\Omega \right) >0$ depend only on listed
parameters. Then the estimate (\ref{2.6}) together with (\ref{2.7}) and (\ref%
{2.8}) is called \emph{pointwise Carleman estimate for the operator }$A_{0}$%
\emph{\ with the CWF }$\varphi _{\lambda }^{2}\left( x\right) $\emph{\ in
the domain }$\Omega _{\alpha }.$

\subsection{Theorems}

\label{sec:3.3}

Let $R>0$ be an arbitrary number. We now specify the ball $B\left( R\right) $
as 
\begin{equation}
B\left( R\right) =\left\{ u\in H_{0}^{k}\left( \Omega \right) :\left\Vert
u\right\Vert _{H^{k}\left( \Omega \right) }<R\right\} .  \label{1}
\end{equation}%
To solve the Cauchy Problem 2, we take into account (\ref{2.1007}) and
consider the following minimization problem:

\textbf{Minimization Problem}. \emph{Assume that the operator }$A_{0}$\emph{%
\ satisfies conditions (\ref{2.6})-(\ref{2.8}). Let }$\beta \in \left(
0,1\right) $\emph{\ be the regularization parameter. Minimize with respect
to the function }$v\in \overline{B}\left( R\right) $ \emph{the functional }$%
J_{\lambda ,\beta }\left( v,F\right) $\emph{, where}%
\begin{equation}
J_{\lambda ,\beta }\left( v,F\right) =e^{-2\lambda \left( \alpha
+\varepsilon \right) }\dint\limits_{\Omega }\left[ A\left( v+F\right) \right]
^{2}\varphi _{\lambda }^{2}dx+\beta \left\Vert v\right\Vert _{H^{k}\left(
\Omega \right) }^{2}.  \label{2.1005}
\end{equation}

The multiplier $e^{-2\lambda \left( \alpha +\varepsilon \right) }$ is
introduced to balance two terms in the right hand side of (\ref{2.1005}).
Below \textquotedblleft the Frech\'{e}t derivative $J_{\lambda ,\beta
}^{\prime }\left( v,F\right) "$ means the Frech\'{e}t derivative of the
functional $J_{\lambda ,\beta }\left( v,F\right) $ with respect to $v$.
Also, below $\left[ ,\right] $ denotes the scalar product in $H^{k}\left(
\Omega \right) .$

\textbf{Theorem 3.1}. \emph{The functional }$J_{\lambda ,\beta }\left(
v,F\right) $\emph{\ has the Frech\'{e}t derivative }$J_{\lambda ,\beta
}^{\prime }\left( v,F\right) \in H_{0}^{k}\left( \Omega \right) $\emph{\ for 
}$v\in B\left( 2R\right) .$\emph{\ This derivative satisfies the Lipschitz
continuity condition}%
\begin{equation}
\left\Vert J_{\lambda ,\beta }^{\prime }\left( v_{1},F\right) -J_{\lambda
,\beta }^{\prime }\left( v_{2},F\right) \right\Vert _{H^{k}\left( \Omega
\right) }\leq L\left\Vert v_{1}-v_{2}\right\Vert _{H^{k}\left( \Omega
\right) },\forall v_{1},v_{2}\in \overline{B}\left( R\right) ,
\label{2.1010}
\end{equation}%
\emph{where the constant }$L=L\left( R,A,F,\Omega ,\lambda ,\alpha
,\varepsilon ,\beta \right) >0$\emph{\ depends only on listed parameters.}

As to Theorem 3.2, we note that since $e^{-\lambda \varepsilon }<<1$ for
sufficiently large $\lambda ,$ then the requirement of this theorem $\beta
\in \left[ e^{-\lambda \varepsilon },1\right) $ enables the regularization
parameter $\beta $ to change from being very small and up to the unity.

\textbf{Theorem 3.2}. \emph{Assume that the operator }$A_{0}$\emph{\ admits
the pointwise\ Carleman estimate (\ref{2.6})-(\ref{2.8}) in the domain }$%
\Omega _{\alpha }$\emph{. Then there exists a sufficiently large number }$%
\lambda _{1}=\lambda _{1}\left( R,A,F,\Omega \right) >\lambda _{0}\left(
A_{0},\Omega \right) >1$\emph{\ and a number }$C_{2}=C_{2}\left(
R,A,F,\Omega \right) >0,$\ \emph{both} \emph{depending only on listed
parameters, such that for all }$\lambda \geq \lambda _{1}$\emph{\ and for
every }$\beta \in \left[ e^{-\lambda \varepsilon },1\right) $\emph{\ the
functional }$J_{\lambda ,\beta }\left( v,F\right) $\emph{\ is strictly
convex on the ball }$\overline{B}\left( R\right) ,$%
\begin{equation}
J_{\lambda ,\beta }\left( v_{2},F\right) -J_{\lambda ,\beta }\left(
v_{1},F\right) -J_{\lambda ,\beta }^{\prime }\left( v_{1},F\right) \left(
v_{2}-v_{1}\right)  \label{2.300}
\end{equation}%
\begin{equation*}
\geq C_{2}e^{2\lambda \varepsilon }\left\Vert v_{2}-v_{1}\right\Vert
_{H^{1}\left( \Omega _{\alpha +2\varepsilon }\right) }^{2}+\frac{\beta }{2}%
\left\Vert v_{2}-v_{1}\right\Vert _{H^{k}\left( \Omega \right) }^{2},\text{ }%
\forall v_{1},v_{2}\in \overline{B}\left( R\right) .
\end{equation*}

To minimize the functional (\ref{2.1005}) on the set $\overline{B}\left(
R\right) $, we apply the gradient projection method. Let $P_{\overline{B}%
\left( R\right) }:H_{0}^{k}\left( \Omega \right) \rightarrow \overline{B}%
\left( R\right) $ be the projection operator of the space $H_{0}^{k}\left(
\Omega \right) $ in the closed ball $\overline{B}\left( R\right) $ (Lemma
2.2). Let an arbitrary function $v_{0}\in \overline{B}\left( R\right) $ be
our starting point for iterations of this method. Let the step size of the
gradient method be $\gamma >0$. Consider the sequence $\left\{ v_{n}\right\}
_{n=0}^{\infty }$, 
\begin{equation}
v_{n+1}=P_{\overline{B}\left( R\right) }\left( v_{n}-\gamma J_{\lambda
,\beta }^{\prime }\left( v_{n},F\right) \right) ,n=0,1,2,...  \label{2.200}
\end{equation}%
For brevity, we do not indicate here the dependence of functions $v_{n}$ on
parameters $\lambda ,\beta ,\gamma $.

\textbf{Theorem 3.3}. \emph{Suppose that all conditions of Theorem 3.2 are
satisfied. Choose a number }$\lambda \geq \lambda _{1}.$\emph{\ Let the
regularization parameter }$\beta \in \left[ e^{-\lambda \varepsilon
},1\right) .$\emph{\ Then there exists a point }$v_{\min }\in \overline{B}%
\left( R\right) $\emph{\ of the relative minimum of the functional }$%
J_{\lambda ,\beta }\left( v\right) $\emph{\ on the set }$\overline{B}\left(
R\right) .$ \emph{Furthermore, }$v_{\min }$\emph{\ is also the unique point
of the global minimum of this functional on }$\overline{B}\left( R\right) .$%
\emph{\ Consider the sequence (\ref{2.200}), where }$v_{0}\in \overline{B}%
\left( R\right) $\emph{\ is an arbitrary point of the closed ball }$%
\overline{B}\left( R\right) $\emph{. Then there exist a sufficiently small
number }$\gamma =\gamma \left( R,A,F,\Omega ,\alpha ,\varepsilon ,\beta
,\lambda \right) \in \left( 0,1\right) $ \emph{and a number }$q\left( \gamma
\right) \in \left( 0,1\right) ,$ \emph{both depending only on listed
parameters,\ such that the sequence (\ref{2.200}) converges to the point }$%
v_{\min },$%
\begin{equation}
\left\Vert v_{n+1}-v_{\min }\right\Vert _{H^{k}\left( \Omega \right) }\leq
q^{n}\left( \gamma \right) \left\Vert v_{0}-v_{\min }\right\Vert
_{H^{k}\left( \Omega \right) },n=0,1,2,...  \label{2.201}
\end{equation}

Following the regularization theory \cite{BK,T}, the next natural question
to address is whether regularized solutions converge to the exact solution
(if it exists) for some values of the parameter $\lambda =\lambda \left(
\delta \right) $ if the level of the error $\delta $ in the Cauchy data $%
g_{0},g_{1}$ tends to zero. Since functions $g_{0},g_{1}$ generate the
function $F$, we consider the error only in $F$. Following \ one of concepts
of the regularization theory, we assume now the existence of the exact
solution $v^{\ast }\in H_{0}^{k}\left( \Omega \right) $\emph{\ }of the
problem (\ref{2.1007}), which\emph{\ }satisfies the following
conditions: 
\begin{equation}
A\left( v^{\ast }+F^{\ast }\right) =0,  \label{2.203}
\end{equation}%
\begin{equation}
v^{\ast }\in B\left( R\right) ,  \label{2.205}
\end{equation}%
where the function $F^{\ast }\in H^{k+1}\left( \Omega \right) $ is generated
by the exact (i.e. noiseless) Cauchy data $g_{0}^{\ast }\left( x\right) $
and $g_{1}^{\ast }\left( x\right) .$ We assume that 
\begin{equation}
\left\Vert F-F^{\ast }\right\Vert _{H^{k+1}\left( \Omega \right) }\leq
\delta ,  \label{2.208}
\end{equation}%
where $\delta \in \left( 0,1\right) $ is a sufficiently small number
characterizing the level of the error in the data. The construction (\ref%
{2.203})-(\ref{2.208}) corresponds well with the regularization theory \cite%
{BK,LRS,T}. First, consider the case when the data are noiseless, i.e. when $%
\delta =0.$

\textbf{Theorem 3.4}. \emph{Suppose that all conditions of Theorem 3.2 are
satisfied. Choose a number }$\lambda ^{\ast }=\lambda ^{\ast }\left(
R,A,F^{\ast },\Omega \right) >\lambda _{0}$ \emph{such that estimate (\ref%
{2.300}) is valid for }$J_{\lambda ,\beta }\left( v,F^{\ast }\right) $\emph{%
\ for all }$\lambda \geq \lambda ^{\ast }.$\emph{\ Let the level of the
error in the data be }$\delta =0$\emph{. Choose }$\lambda \geq \lambda
^{\ast }$ \emph{and }$\beta =e^{-\lambda \varepsilon }.$\emph{\ Let }$%
v_{\min }\in \overline{B}\left( R\right) $\emph{\ be the point of the unique
global minimum on }$\overline{B}\left( R\right) $\emph{\ of the functional }$%
J_{\lambda ,\beta }\left( v,F^{\ast }\right) $\emph{\ (Theorem 3.3). Then
there exists a constant }$C_{3}=C_{3}\left( R,A,F^{\ast },\Omega \right) >0$%
\emph{\ depending only on listed parameters such that }%
\begin{equation}
\left\Vert v^{\ast }-v_{\min }\right\Vert _{H^{1}\left( \Omega _{\alpha
+2\varepsilon }\right) }\leq C_{3}\exp \left( -3\lambda \varepsilon
/2\right) .  \label{2.209}
\end{equation}%
\emph{Furthermore, let }$\left\{ v_{n}\right\} _{n=0}^{\infty }$\emph{\ be
the sequence (\ref{2.200}) where the number }

$\gamma =\gamma \left( R,A,F^{\ast },\Omega ,\alpha ,\varepsilon ,\beta
,\lambda \right) \in \left( 0,1\right) $\emph{\ is the same as in Theorem
3.3. Then with the same constant }$q\left( \gamma \right) \in \left(
0,1\right) $\emph{\ as in Theorem 3.3 the following estimate holds:}%
\begin{equation}
\left\Vert v^{\ast }-v_{n+1}\right\Vert _{H^{1}\left( \Omega _{\alpha
+2\varepsilon }\right) }\leq C_{3}\exp \left( -3\lambda \varepsilon
/2\right) +q^{n}\left( \gamma \right) \left\Vert v_{0}-v_{\min }\right\Vert
_{H^{k}\left( \Omega \right) },n=0,1,2,...  \label{2.2009}
\end{equation}

Let $m$ be the number in (\ref{2.30}). Denote%
\begin{equation}
\theta =\min \left( \frac{\varepsilon }{4m},\frac{1}{2}\right) .  \label{100}
\end{equation}%
Theorem 3.5 estimates the rate of convergence of minimizers $v_{\min }$ to
the exact solution $v^{\ast }$ in the norm of the space $H^{1}\left( \Omega
_{\alpha +2\varepsilon }\right) .$

\textbf{Theorem 3.5. }\emph{Let all conditions of Theorem 3.2 hold. Let the
number }$\lambda _{1}=\lambda _{1}\left( R,A,F,\Omega \right) >\lambda _{0}$%
\emph{\ be the same as in Theorem 3.2 and let }$\theta $ \emph{be the number
defined in (\ref{100}). Let the number }$\delta _{0}\in \left( 0,1\right) $%
\emph{\ be so small that }$\delta _{0}^{-1/\left( 2m\right) }>e^{\lambda
_{1}}.$\emph{\ Let }$\delta \in \left( 0,\delta _{0}\right) $\emph{\ be the
level of the error in the function }$F,$\emph{\ i.e. let (\ref{2.208}) be
valid. Choose }$\lambda =\lambda \left( \delta \right) =\ln \left( \delta
^{-1/\left( 2m\right) }\right) >\lambda _{1}$ \emph{and} $\beta =e^{-\lambda
\left( \delta \right) \varepsilon }.$\emph{\ Let }$v_{\min }\in \overline{B}%
\left( R\right) $\emph{\ be the point of the unique global minimum on }$%
\overline{B}\left( R\right) $\emph{\ of the functional }$J_{\lambda ,\beta
}\left( v,F\right) $\emph{\ (Theorem 3.3). Then there exists a constant }$%
C_{4}=C_{4}\left( R,A,F,\Omega \right) >0$\emph{\ depending only on listed
parameters such that}%
\begin{equation}
\left\Vert v^{\ast }-v_{\min }\right\Vert _{H^{1}\left( \Omega _{\alpha
+2\varepsilon }\right) }\leq C_{4}\delta ^{\theta }.  \label{2}
\end{equation}%
\emph{Next, let let }$\left\{ v_{n}\right\} _{n=0}^{\infty }$\emph{\ be the
sequence (\ref{2.200}), where the number }

$\gamma =\gamma \left( R,A,F,\Omega ,\alpha ,\varepsilon ,\beta ,\delta
\right) \in \left( 0,1\right) $\emph{\ is the same as in Theorem 3.3. Then
with the same constant }$q\left( \gamma \right) \in \left( 0,1\right) $\emph{%
\ as in Theorem 3.3 the following estimate holds:}%
\begin{equation}
\left\Vert v^{\ast }-v_{n+1}\right\Vert _{H^{1}\left( \Omega _{\alpha
+2\varepsilon }\right) }\leq C_{4}\delta ^{\theta }+q^{n}\left( \gamma
\right) \left\Vert v_{0}-v_{\min }\right\Vert _{H^{k}\left( \Omega \right)
},n=0,1,2,...  \label{3}
\end{equation}

\textbf{Remarks 3.1}:

\begin{enumerate}
\item We point out that, compared with previous publications \cite%
{BKconv,Klib97,Kpar,KNT,Kl1,KK,KKLY} on the topic of this paper, a
significantly new element of Theorems 3.3-3.5 is that now the existence of
the global minimum $v_{\min }$ is asserted rather than assumed. This became
possible because of results of convex analysis of section 2.

\item Even though we estimate in (\ref{2.209})-(\ref{2.2009}) only norms in $%
H^{1}\left( \Omega _{\alpha +2\varepsilon }\right) ,$ this seems to be
sufficient for computations, see section 9. It follows from the combination
of Theorems 3.2-3.5 that the optimization procedure (\ref{2.200}) represents
a globally convergent numerical method for the Cauchy Problem 2. Here the
global convergence is understood as described in section 1.

\item Theorem 3.3 follows immediately from Theorems 2.1, 3.1 and 3.2. Hence,
we do not prove Theorem 3.3 here. However, we still need to prove all other
theorems, since their proofs are essentially different from proofs of
similar theorems in \cite{Kl1}. These differences are caused by two factors.
First, we now introduce the function $F$ in (\ref{2.1005}), which was not
the case of previous publications. Second, we now integrate in the first
term in the right hand side of (\ref{2.1005}) over the entire domain $\Omega
.$ On the other hand, the integration was carried out over the subdomain $%
\Omega _{\alpha }$ in \cite{Kl1}.
\end{enumerate}

\section{Proof of Theorem 3.1}

\label{sec:4}

In this proof $L=L\left( R,A,F,\Omega ,\alpha ,\varepsilon ,\beta ,\lambda
\right) >0$ denotes different numbers depending only on listed parameters.
Let $v_{1},v_{2}\in B\left( 2R\right) $ be two arbitrary functions. Denote $%
h=v_{2}-v_{1}.$ Hence, $h\in H_{0}^{k}\left( \Omega \right) .$ Let 
\begin{equation}
D=\left( A\left( v_{2}+F\right) \right) ^{2}-\left( A\left( v_{1}+F\right)
\right) ^{2}.  \label{10.1}
\end{equation}%
By the Lagrange formula 
\begin{equation}
f\left( y+z\right) =f\left( y\right) +f^{\prime }\left( y\right) z+\frac{%
z^{2}}{2}f^{\prime \prime }\left( \eta \right) ,\forall y,z\in \mathbb{R}%
,\forall f\in C^{2}\left( \mathbb{R}\right) ,  \label{10.2}
\end{equation}%
where $\eta =\eta \left( y,z\right) $ is a number located between numbers $y$
and $y+z$. By (\ref{2.10005})%
\begin{equation}
\left\Vert h\right\Vert _{C^{1}\left( \overline{\Omega }\right) }=\left\Vert
v_{2}-v_{1}\right\Vert _{C^{1}\left( \overline{\Omega }\right) }\leq 4CR.
\label{10.3}
\end{equation}%
Hence, using (\ref{2.100})-(\ref{2.1002}), (\ref{10.2}) and (\ref{10.3}), we
obtain 
\begin{equation*}
A_{1}\left( x,\nabla \left( v_{2}+F\right) ,v_{2}+F\right) =A_{1}\left(
x,\nabla \left( v_{1}+F+h\right) ,v_{1}+F+h\right)
\end{equation*}%
\begin{equation*}
=A_{1}\left( x,\nabla v_{1}+\nabla F,v_{1}+F\right) +
\end{equation*}%
\begin{equation*}
+\dsum\limits_{i=1}^{n}\partial _{v_{x_{i}}}A_{1}\left( x,\nabla
v_{1}+\nabla F,v_{1}+F\right) h_{x_{i}}+\partial _{v}A_{1}\left( x,\nabla
v_{1}+\nabla F,v_{1}+F\right) h
\end{equation*}%
\begin{equation*}
+P\left( x,\nabla v_{1}+\nabla F,v_{1}+F,\nabla h,h\right) ,
\end{equation*}%
where the function $P$ satisfies the following estimate%
\begin{equation}
\left\vert P\left( x,\nabla v_{1}+\nabla F,v_{1}+F,\nabla h,h\right)
\right\vert \leq K\left( \left( \nabla h\right) ^{2}+h^{2}\right) ,\forall
x\in \overline{\Omega },\forall v_{1}\in B\left( 2R\right) ,  \label{10.4}
\end{equation}%
where the constant $K=K\left( R,F,\Omega \right) >0$ depends only on listed
parameters. Hence, 
\begin{equation*}
A\left( v_{2}+F\right) =A_{0}\left( v_{1}+F+h\right) +A_{1}\left( x,\nabla
\left( v_{1}+F+h\right) ,v_{1}+F+h\right) =A\left( v_{1}+F\right)
\end{equation*}%
\begin{equation*}
+\left[ A_{0}\left( h\right) +\dsum\limits_{i=1}^{n}\partial
_{v_{x_{i}}}A_{1}\left( x,\nabla v_{1}+\nabla F,v_{1}+F\right)
h_{x_{i}}+\partial _{v}A_{1}\left( x,\nabla v_{1}+\nabla F,v_{1}+F\right) h%
\right]
\end{equation*}%
\begin{equation*}
+P\left( x,\nabla u_{1},u_{1},h\right) .
\end{equation*}%
Hence, by (\ref{10.1})%
\begin{equation*}
D=2A\left( v_{1}+F\right) \times
\end{equation*}%
\begin{equation}
\left[ A_{0}\left( h\right) +\dsum\limits_{i=1}^{n}\partial
_{v_{x_{i}}}A_{1}\left( x,\nabla v_{1}+\nabla F,v_{1}+F\right)
h_{x_{i}}+\partial _{v}A_{1}\left( x,\nabla v_{1}+\nabla F,v_{1}+F\right) h%
\right]  \label{10.5}
\end{equation}%
\begin{equation*}
+\left[ A_{0}\left( h\right) +\dsum\limits_{i=1}^{n}\partial
_{v_{x_{i}}}A_{1}\left( x,\nabla v_{1}+\nabla F,v_{1}+F\right)
h_{x_{i}}+\partial _{v}A_{1}\left( x,\nabla v_{1}+\nabla F,v_{1}+F\right) h%
\right] ^{2}
\end{equation*}%
\begin{equation*}
+P^{2}.
\end{equation*}

The expression in the first two lines of (\ref{10.5}) is linear with respect
to $h$. We denote this expression as $Q\left( v_{1}+F\right) \left( h\right)
,$%
\begin{equation}
Q\left( v_{1}+F\right) \left( h\right) =2A\left( v_{1}+F\right) \times
\label{10.50}
\end{equation}%
\begin{equation*}
\left[ A_{0}\left( h\right) +\dsum\limits_{i=1}^{n}\partial
_{v_{x_{i}}}A_{1}\left( x,\nabla v_{1}+\nabla F,v_{1}+F\right)
h_{x_{i}}+\partial _{v}A_{1}\left( x,\nabla v_{1}+\nabla F,v_{1}+F\right) h%
\right] .
\end{equation*}%
Consider the linear functional acting on functions $h\in H_{0}^{k}\left(
\Omega \right) $ as 
\begin{equation}
\widetilde{J}\left( v_{1},F\right) \left( h\right) =\dint\limits_{\Omega
}Q\left( v_{1}+F\right) \left( h\right) \varphi _{\lambda }^{2}dx+2\beta %
\left[ v_{1},h\right] .  \label{10.6}
\end{equation}%
Clearly, $\widetilde{J}\left( v_{1},F\right) \left( h\right)
:H_{0}^{k}\left( \Omega \right) \rightarrow \mathbb{R}$ is a bounded linear
functional. Hence, by the Riesz theorem, there exists a single element $%
M\left( v_{1}\right) \in H_{0}^{k}\left( \Omega \right) $ such that 
\begin{equation}
\widetilde{J}\left( v_{1},F\right) \left( h\right) =\left[ M\left(
v_{1},F\right) ,h\right] ,\forall h\in H_{0}^{k}\left( \Omega \right) .
\label{10.7}
\end{equation}%
Furthermore, 
\begin{equation}
\left\Vert M\left( v_{1},F\right) \right\Vert _{H^{k}\left( \Omega \right)
}=\left\Vert \widetilde{J}\left( v_{1},F\right) \right\Vert .  \label{10.8}
\end{equation}%
Next, since by (\ref{2.10005}) $\left\Vert h\right\Vert _{C^{1}\left( 
\overline{\Omega }\right) }\leq C\left\Vert h\right\Vert _{H^{k}\left(
\Omega \right) },$ then (\ref{2.1005}), (\ref{10.1}), (\ref{10.5}) and (\ref%
{10.6}) imply that 
\begin{equation}
J_{\lambda ,\beta }\left( v_{1}+h,F\right) -J_{\lambda ,\beta }\left(
v_{1},F\right) -\widetilde{J}\left( v_{1},F\right) \left( h\right) =O\left(
\left\Vert h\right\Vert _{H^{k}\left( \Omega \right) }^{2}\right) ,\text{ }
\label{10.9}
\end{equation}%
as $\left\Vert h\right\Vert _{H^{k}\left( \Omega \right) }\rightarrow 0.$
The existence of the Frech\'{e}t derivative $J_{\lambda ,\beta }^{^{\prime
}}\left( v_{1}\right) $ follows from (\ref{10.50})-(\ref{10.9}). Also, for
all $h\in H_{0}^{k}\left( \Omega \right) $ and all $v\in B\left( 2R\right) $ 
\begin{equation}
J_{\lambda ,\beta }^{^{\prime }}\left( v,F\right) \left( h\right) =%
\widetilde{J}\left( v,F\right) \left( h\right) =\dint\limits_{\Omega
}Q\left( v+F\right) \left( h\right) \varphi _{\lambda }^{2}dx+2\beta \left[
v,h\right] ,  \label{10.10}
\end{equation}%
\begin{equation}
J_{\lambda ,\beta }^{^{\prime }}\left( v,F\right) =M\left( v,F\right) \in
H_{0}^{k}\left( \Omega \right) .  \label{10.11}
\end{equation}

We now prove the Lipschitz continuity of the Frech\'{e}t derivative $%
J_{\lambda ,\beta }^{^{\prime }}\left( v,F\right) .$ By (\ref{10.5}), (\ref%
{10.50}), (\ref{10.6}), (\ref{10.10}) and (\ref{10.11}) we should analyze
the following expression for all $v_{1},v_{2}\in \overline{B}\left( R\right) 
$ and for all $h\in H_{0}^{k}\left( \Omega \right) :$%
\begin{equation*}
Y\left( v_{1},h\right) -Y\left( v_{2},h\right) =
\end{equation*}%
\begin{equation*}
2A\left( v_{1}+F\right) \times
\end{equation*}%
\begin{equation}
\left[ A_{0}\left( h\right) +\dsum\limits_{i=1}^{n}\partial
_{v_{x_{i}}}A_{1}\left( x,\nabla v_{1}+\nabla F,v_{1}+F\right)
h_{x_{i}}+\partial _{v}A_{1}\left( x,\nabla v_{1}+\nabla F,v_{1}+F\right) h%
\right]  \label{10.12}
\end{equation}%
\begin{equation*}
-2A\left( v_{2}+F\right) \times
\end{equation*}%
\begin{equation*}
\left[ A_{0}\left( h\right) +\dsum\limits_{i=1}^{n}\partial
_{v_{x_{i}}}A_{1}\left( x,\nabla v_{2}+\nabla F,v_{2}+F\right)
h_{x_{i}}+\partial _{v}A_{1}\left( x,\nabla v_{2}+\nabla F,v_{2}+F\right) h%
\right] .
\end{equation*}%
We have 
\begin{equation*}
Y\left( v_{1},h\right) -Y\left( v_{2},h\right) =2\left( A\left(
v_{1}+F\right) -A\left( v_{2}+F\right) \right) \times
\end{equation*}%
\begin{equation*}
\left[ A_{0}\left( h\right) +\dsum\limits_{i=1}^{n}\partial
_{v_{x_{i}}}A_{1}\left( x,\nabla v_{1}+\nabla F,v_{1}+F\right)
h_{x_{i}}+\partial _{v}A_{1}\left( x,\nabla v_{1}+\nabla F,v_{1}+F\right) h%
\right]
\end{equation*}%
\begin{equation}
+2A\left( v_{2}+F\right) \times  \label{10.13}
\end{equation}%
\begin{equation*}
\left[ \dsum\limits_{i=1}^{n}\partial _{v_{x_{i}}}A_{1}\left( x,\nabla
v_{1}+\nabla F,v_{1}+F\right) h_{x_{i}}+\partial _{v}A_{1}\left( x,\nabla
v_{1}+\nabla F,v_{1}+F\right) h\right]
\end{equation*}%
\begin{equation*}
-2A\left( v_{2}+F\right) \times
\end{equation*}%
\begin{equation*}
\left[ \dsum\limits_{i=1}^{n}\partial _{v_{x_{i}}}A_{1}\left( x,\nabla
v_{2}+\nabla F,v_{2}+F\right) h_{x_{i}}+\partial _{v}A_{1}\left( x,\nabla
v_{2}+\nabla F,v_{2}+F\right) h\right] .
\end{equation*}

First, using (\ref{2.100}) and (\ref{10.2}), we obtain 
\begin{equation}
2\left( A\left( v_{1}+F\right) -A\left( v_{2}+F\right) \right) =2A_{0}\left(
v_{1}-v_{2}\right)  \label{10.14}
\end{equation}%
\begin{equation*}
+2\left[ \dsum\limits_{i=1}^{n}\partial _{v_{x_{i}}}A_{1}\left( x,\nabla
v_{2}+\nabla F,v_{2}+F\right) \left( v_{1}-v_{2}\right) _{x_{i}}+\partial
_{v}A_{1}\left( x,\nabla v_{2}+\nabla F,v_{2}+F\right) \left(
v_{1}-v_{2}\right) \right]
\end{equation*}%
\begin{equation*}
+Y_{1}\left( x,v_{1},v_{2}\right) ,
\end{equation*}%
where 
\begin{equation}
\left\vert Y_{1}\left( x,v_{1},v_{2}\right) \right\vert \leq L\left\Vert
v_{1}-v_{2}\right\Vert _{C^{1}\left( \overline{\Omega }\right) }^{2}\leq
L\left\Vert v_{1}-v_{2}\right\Vert _{H^{k}\left( \Omega \right) },\forall
v_{1},v_{2}\in \overline{B}\left( R\right) .  \label{10.15}
\end{equation}%
Thus, (\ref{10.14}) and (\ref{10.15}) imply that the modulus of the
expression in the first two lines of (\ref{10.13}) can be estimated from the
above via $Y_{2},$ where 
\begin{equation}
Y_{2}\leq L\left\Vert v_{1}-v_{2}\right\Vert _{H^{k}\left( \Omega \right)
}\left\Vert h\right\Vert _{H^{k}\left( \Omega \right) },\forall
v_{1},v_{2}\in \overline{B}\left( R\right) ,\forall h\in H_{0}^{k}\left(
\Omega \right) .  \label{10.16}
\end{equation}

Estimate now from the above the modulus of the expression in the lines
number 3-6 of (\ref{10.13}). By (\ref{10.2})%
\begin{equation*}
\partial _{v}A_{1}\left( x,\nabla v_{1}+\nabla F,v_{1}+F\right) h-\partial
_{v}A_{1}\left( x,\nabla v_{2}+\nabla F,v_{2}+F\right) h
\end{equation*}%
\begin{equation*}
=\partial _{v}^{2}A_{1}\left( x,\nabla v_{2}+\nabla F,v_{2}+F\right) h\left(
v_{1}-v_{2}\right) +\frac{\left( v_{1}-v_{2}\right) ^{2}}{2}h\partial
_{v}^{3}A_{1}\left( x,\nabla v_{2}+\nabla F,\xi \left( x\right) +F\right) ,
\end{equation*}%
where the point $\xi \left( x\right) $ is located between points $%
v_{1}\left( x\right) $ and $v_{2}\left( x\right) .$ Similar formulas are
valid of course for terms%
\begin{equation*}
\dsum\limits_{i=1}^{n}\partial _{v_{x_{i}}}A_{1}\left( x,\nabla v_{1}+\nabla
F,v_{1}+F\right) h_{x_{i}}-\dsum\limits_{i=1}^{n}\partial
_{v_{x_{i}}}A_{1}\left( x,\nabla v_{2}+\nabla F,v_{2}+F\right) h_{x_{i}}.
\end{equation*}%
Hence, the modulus of the expression in lines number 3-6 of (\ref{10.13})
can be estimated from the above similarly with (\ref{10.16}) via $Y_{3},$
where%
\begin{equation}
Y_{3}\leq L\left\Vert v_{1}-v_{2}\right\Vert _{H^{k}\left( \Omega \right)
}\left\Vert h\right\Vert _{H^{k}\left( \Omega \right) },\forall
v_{1},v_{2}\in \overline{B}\left( R\right) ,\forall h\in H_{0}^{k}\left(
\Omega \right) .  \label{10.17}
\end{equation}%
Thus, (\ref{10.50}) and (\ref{10.10})-(\ref{10.17}) imply that 
\begin{equation*}
\left\vert J_{\lambda ,\beta }^{^{\prime }}\left( v_{1},F\right) \left(
h\right) -J_{\lambda ,\beta }^{^{\prime }}\left( v_{2},F\right) \left(
h\right) \right\vert \leq L\left\Vert v_{1}-v_{2}\right\Vert _{H^{k}\left(
\Omega \right) }\left\Vert h\right\Vert _{H^{k}\left( \Omega \right) }.
\end{equation*}%
for all $v_{1},v_{2}\in \overline{B}\left( R\right) $ and for all $h\in
H_{0}^{k}\left( \Omega \right) .$ This, in turn implies (\ref{2.1010}). $%
\square $

\section{Proof of Theorem 3.2}

\label{sec:5}

In this proof $C_{2}=C_{2}\left( R,A,F,\Omega \right) >C_{1}>0$ denotes
different constants depending only on listed parameters. Here $%
C_{1}=C_{1}\left( A_{0},\Omega \right) >0$ is the constant of the pointwise
Carleman estimate (\ref{2.6})-(\ref{2.8}). For two arbitrary points $%
v_{1},v_{2}\in \overline{B}\left( R\right) $ let again $h=v_{2}-v_{1}$ and
let $D$ be the same as in (\ref{10.1}). Denote $S=D-Q\left( v_{1}+F\right)
\left( h\right) ,$ where $Q\left( v_{1}+F\right) \left( h\right) $ is given
in (\ref{10.50}) and it is linear, with respect to $h.$ Then, using (\ref%
{10.4})-(\ref{10.50}) and the Cauchy-Schwarz inequality, we obtain%
\begin{equation*}
S\geq \frac{1}{2}\left( A_{0}h\right) ^{2}-C_{2}\left( \left( \nabla
h\right) ^{2}+h^{2}\right) ,\forall x\in \Omega .
\end{equation*}%
Hence, using (\ref{10.9}) and (\ref{10.10}), we obtain%
\begin{equation*}
J_{\lambda ,\beta }\left( v_{1}+h,F\right) -J_{\lambda ,\beta }\left(
v_{1},F\right) -J_{\lambda ,\beta }^{\prime }\left( v_{1},F\right) \left(
h\right)
\end{equation*}%
\begin{equation}
\geq \frac{1}{2}e^{-2\lambda \left( \alpha +\varepsilon \right)
}\dint\limits_{\Omega }\left( A_{0}h\right) ^{2}\varphi _{\lambda
}^{2}dx-C_{2}e^{-2\lambda \left( \alpha +\varepsilon \right)
}\dint\limits_{\Omega }\left( \left( \nabla h\right) ^{2}+h^{2}\right)
\varphi _{\lambda }^{2}dx+\beta \left\Vert h\right\Vert _{H^{k}\left( \Omega
\right) }^{2}.  \label{11.1}
\end{equation}%
Since $\Omega _{\alpha }\subset \Omega ,$ then 
\begin{equation}
e^{-2\lambda \left( \alpha +\varepsilon \right) }\dint\limits_{\Omega
}\left( A_{0}h\right) ^{2}\varphi _{\lambda }^{2}dx\geq e^{-2\lambda \left(
\alpha +\varepsilon \right) }\dint\limits_{\Omega _{\alpha }}\left(
A_{0}h\right) ^{2}\varphi _{\lambda }^{2}dx.  \label{11.2}
\end{equation}%
Next, 
\begin{equation}
-C_{2}e^{-2\lambda \left( \alpha +\varepsilon \right) }\dint\limits_{\Omega
}\left( \left( \nabla h\right) ^{2}+h^{2}\right) \varphi _{\lambda
}^{2}dx=-C_{2}e^{-2\lambda \left( \alpha +\varepsilon \right)
}\dint\limits_{\Omega _{\alpha }}\left( \left( \nabla h\right)
^{2}+h^{2}\right) \varphi _{\lambda }^{2}dx  \label{11.3}
\end{equation}%
\begin{equation*}
-C_{2}e^{-2\lambda \left( \alpha +\varepsilon \right) }\dint\limits_{\Omega
\diagdown \Omega _{\alpha }}\left( \left( \nabla h\right) ^{2}+h^{2}\right)
\varphi _{\lambda }^{2}dx.
\end{equation*}%
Since by (\ref{2.0}) and (\ref{2.2}) $\varphi _{\lambda }^{2}\left( x\right)
<\exp \left( 2\lambda \alpha \right) $ for $x\in \Omega \diagdown \Omega
_{\alpha },$ then%
\begin{equation}
-C_{2}e^{-2\lambda \left( \alpha +\varepsilon \right) }\dint\limits_{\Omega
\diagdown \Omega _{\alpha }}\left( \left( \nabla h\right) ^{2}+h^{2}\right)
\varphi _{\lambda }^{2}dx\geq -C_{2}e^{-2\lambda \varepsilon
}\dint\limits_{\Omega \diagdown \Omega _{\alpha }}\left( \left( \nabla
h\right) ^{2}+h^{2}\right) dx.  \label{11.4}
\end{equation}

Integrate (\ref{2.6}) over the domain $\Omega _{\alpha },$ using the Gauss'
formula, (\ref{2.7}) and (\ref{2.8}). Next, replace $u$ with $h$ in the
resulting formula. Even though there is no guarantee that $h\in C^{2}\left( 
\overline{\Omega }_{\alpha }\right) ,$ still density arguments ensure that
the resulting inequality remains true. Hence, taking into account (\ref{2.0}%
)-(\ref{2.3}), (\ref{11.1}) and (\ref{11.2}), we obtain 
\begin{equation}
\frac{1}{2}e^{-2\lambda \left( \alpha +\varepsilon \right)
}\dint\limits_{\Omega _{\alpha }}\left( A_{0}h\right) ^{2}\varphi _{\lambda
}^{2}dx\geq C_{2}e^{-2\lambda \left( \alpha +\varepsilon \right)
}\dint\limits_{\Omega _{\alpha }}\left( \lambda \left( \nabla h\right)
^{2}+\lambda ^{3}h^{2}\right) \varphi _{\lambda }^{2}dx  \label{11.5}
\end{equation}%
\begin{equation*}
-C_{2}\lambda ^{3}e^{-2\lambda \varepsilon }\dint\limits_{\psi _{\alpha
}}\left( \left( \nabla h\right) ^{2}+h^{2}\right) dS,\forall \lambda \geq
\lambda _{0}.
\end{equation*}%
Since $k\geq 2,$ then the trace theorem implies that 
\begin{equation}
C_{2}\lambda ^{3}e^{-2\lambda \varepsilon }\dint\limits_{\psi _{\alpha
}}\left( \left( \nabla h\right) ^{2}+h^{2}\right) dx\leq C_{2}\lambda
^{3}e^{-2\lambda \varepsilon }\left\Vert h\right\Vert _{H^{k}\left( \Omega
\right) }^{2}.  \label{11.6}
\end{equation}%
Also,%
\begin{equation}
C_{2}e^{-2\lambda \varepsilon }\dint\limits_{\Omega \diagdown \Omega
_{\alpha }}\left( \left( \nabla h\right) ^{2}+h^{2}\right) dx\leq
C_{2}e^{-2\lambda \varepsilon }\left\Vert h\right\Vert _{H^{k}\left( \Omega
\right) }^{2}.  \label{11.7}
\end{equation}%
Since $\beta \geq e^{-\lambda \varepsilon },$ then (\ref{11.6}) and (\ref%
{11.7}) imply that for sufficiently large

$\lambda _{1}=\lambda _{1}\left( R,A,F,\Omega ,\alpha ,\varepsilon ,\beta
\right) >\lambda _{0}$ and for $\lambda \geq \lambda _{1}$%
\begin{equation}
-C_{2}\lambda ^{3}e^{-2\lambda \varepsilon }\dint\limits_{\psi _{\alpha
}}\left( \left( \nabla h\right) ^{2}+h^{2}\right) dx-C_{2}e^{-2\lambda
\varepsilon }\dint\limits_{\Omega \diagdown \Omega _{\alpha }}\left( \left(
\nabla h\right) ^{2}+h^{2}\right) dx\geq -\frac{\beta }{2}\left\Vert
h\right\Vert _{H^{k}\left( \Omega \right) }^{2}.  \label{11.8}
\end{equation}%
Also, for $\lambda \geq \lambda _{1}$%
\begin{equation}
C_{2}e^{-2\lambda \left( \alpha +\varepsilon \right) }\dint\limits_{\Omega
_{\alpha }}\left( \lambda \left( \nabla h\right) ^{2}+\lambda
^{3}h^{2}\right) \varphi _{\lambda }^{2}dx-C_{2}e^{-2\lambda \left( \alpha
+\varepsilon \right) }\dint\limits_{\Omega _{\alpha }}\left( \left( \nabla
h\right) ^{2}+h^{2}\right) \varphi _{\lambda }^{2}dx  \label{11.9}
\end{equation}%
\begin{equation*}
\geq \frac{1}{2}C_{2}e^{-2\lambda \left( \alpha +\varepsilon \right)
}\dint\limits_{\Omega _{\alpha }}\left( \lambda \left( \nabla h\right)
^{2}+\lambda ^{3}h^{2}\right) \varphi _{\lambda }^{2}dx.
\end{equation*}%
Hence, using (\ref{11.1}), (\ref{11.3}), (\ref{11.5}), (\ref{11.8}) and (\ref%
{11.9}), we obtain for $\lambda \geq \lambda _{1}$ with a new constant $%
C_{2} $ 
\begin{equation*}
J_{\lambda ,\beta }\left( v_{1}+h,F\right) -J_{\lambda ,\beta }\left(
v_{1},F\right) -J_{\lambda ,\beta }^{\prime }\left( v_{1},F\right) \left(
h\right)
\end{equation*}%
\begin{equation}
\geq C_{2}e^{-2\lambda \left( \alpha +\varepsilon \right)
}\dint\limits_{\Omega _{\alpha }}\left( \lambda \left( \nabla h\right)
^{2}+\lambda ^{3}h^{2}\right) \varphi _{\lambda }^{2}dx+\frac{\beta }{2}%
\left\Vert h\right\Vert _{H^{k}\left( G_{c}\right) }^{2}.  \label{11.10}
\end{equation}%
Next, since $\Omega _{\alpha +2\varepsilon }\subset \Omega _{\alpha }$ and $%
\varphi _{\lambda }^{2}\left( x\right) >e^{2\lambda \left( \alpha
+2\varepsilon \right) }$ for $x\in \Omega _{\alpha +2\varepsilon },$ then (%
\ref{11.10}) implies that for all $v_{1},v_{2}=v_{1}+h\in \overline{B}\left(
R\right) $ 
\begin{equation*}
J_{\lambda ,\beta }\left( v_{1}+h,F\right) -J_{\lambda ,\beta }\left(
v_{1},F\right) -J_{\lambda ,\beta }^{\prime }\left( v_{1},F\right) \left(
h\right) \geq C_{2}e^{2\lambda \varepsilon }\left\Vert h\right\Vert
_{H^{1}\left( \Omega _{\alpha +2\varepsilon }\right) }^{2}+\frac{\beta }{2}%
\left\Vert h\right\Vert _{H^{k}\left( \Omega \right) }^{2}.\text{ \ \ \ \ }%
\square
\end{equation*}

\section{Proof of Theorem 3.4}

\label{sec:6}

Recall that $\lambda \geq \lambda ^{\ast }.$ The existence and uniqueness of
the point $v_{\min }\in \overline{B}\left( R\right) $ of the global minimum
of the functional $J_{\lambda ,\beta }\left( v,F^{\ast }\right) $ follows
immediately from Theorems 2.1, 3.2 and 3.3. Since by (\ref{2.203}) $A\left(
v^{\ast }+F^{\ast }\right) =0$ and by (\ref{2.205}) $v^{\ast }\in B\left(
R\right) ,$ then, using (\ref{2.1005}), we obtain%
\begin{equation}
J_{\lambda ,\beta }\left( v^{\ast },F^{\ast }\right) =\beta \left\Vert
v^{\ast }\right\Vert _{H^{k}\left( \Omega \right) }^{2}\leq \beta R^{2}.
\label{12.1}
\end{equation}%
Next, by (\ref{8.3})%
\begin{equation}
-J_{\lambda ,\beta }^{\prime }\left( v_{\min },F^{\ast }\right) \left(
v^{\ast }-v_{\min }\right) \leq 0.  \label{12.2}
\end{equation}%
Hence, combining (\ref{12.1}) and (\ref{12.2}), we obtain 
\begin{equation}
J_{\lambda ,\beta }\left( v^{\ast },F^{\ast }\right) -J_{\lambda ,\beta
}\left( v_{\min },F^{\ast }\right) -J_{\lambda ,\beta }^{\prime }\left(
v_{\min },F^{\ast }\right) \left( v^{\ast }-v_{\min }\right) \leq \beta
R^{2}.  \label{12.3}
\end{equation}%
Next, combining (\ref{12.3}) with Theorem 3.2 and setting \ $\beta
=e^{-\lambda \varepsilon }$, we obtain (\ref{2.209}). Next, since 
\begin{equation*}
\left\Vert v^{\ast }-v_{n+1}\right\Vert _{H^{1}\left( \Omega _{\alpha
+2\varepsilon }\right) }\leq \left\Vert v^{\ast }-v_{\min }\right\Vert
_{H^{1}\left( \Omega _{\alpha +2\varepsilon }\right) }+\left\Vert v_{\min
}-v_{n+1}\right\Vert _{H^{1}\left( \Omega _{\alpha +2\varepsilon }\right) }
\end{equation*}%
and $\left\Vert v_{\min }-v_{n+1}\right\Vert _{H^{1}\left( \Omega _{\alpha
+2\varepsilon }\right) }\leq \left\Vert v_{\min }-v_{n+1}\right\Vert
_{H^{k}\left( \Omega \right) },$ then (\ref{2.2009}) follows from (\ref%
{2.209}) and (\ref{2.201}). $\square $

\section{Proof of Theorem 3.5}

\label{sec:7}

In this proof $C_{4}=C_{4}\left( R,A,F,\Omega \right) >0$ denotes different
constants depending only on listed parameters. Since functions $F,F^{\ast
}\in H^{k+1}\left( \Omega \right) ,$ then, as it was noticed in subsection
3.1, 
\begin{equation}
A\left( F\right) ,A\left( F^{\ast }\right) \in C\left( \overline{\Omega }%
\right) .  \label{7.1}
\end{equation}%
It follows from (\ref{2.100}), (\ref{2.203})-(\ref{2.208}), (\ref{10.2}) and
(\ref{7.1}) that 
\begin{equation*}
A\left( v^{\ast }+F\right) =A\left( v^{\ast }+F^{\ast }+\left( F-F^{\ast
}\right) \right) =A\left( v^{\ast }+F^{\ast }\right) +\widetilde{A}\left(
v^{\ast },F-F^{\ast }\right) =\widetilde{A}\left( v^{\ast },F-F^{\ast
}\right) ,
\end{equation*}%
where $\left\vert \widetilde{A}\left( v^{\ast },F-F^{\ast }\right)
\right\vert \leq C_{4}\delta ,\forall x\in \overline{\Omega }.$ Hence,
recalling that $v^{\ast }\in B\left( R\right) $ and applying (\ref{2.4}) and
(\ref{2.1005}), we obtain%
\begin{equation}
J_{\lambda ,\beta }\left( v^{\ast },F\right) \leq C_{4}\left( \delta
^{2}e^{2\lambda m}+\beta \right) .  \label{7.3}
\end{equation}%
Recall that $\lambda \geq \lambda _{1}$. Let $v_{\min }\in \overline{B}%
\left( R\right) $ be the unique point of the global minimum of the
functional $J_{\lambda ,\beta }\left( v,F\right) $ on the set $\overline{B}%
\left( R\right) .$ Since (\ref{12.2}) is valid, then, using Theorem 3.2 and (%
\ref{7.3}), we obtain%
\begin{equation}
\left\Vert v^{\ast }-v_{\min }\right\Vert _{H^{1}\left( \Omega _{\alpha
+2\varepsilon }\right) }\leq C_{4}\left( \delta e^{\lambda m}+\sqrt{\beta }%
\right) .  \label{7.4}
\end{equation}%
Choose $\lambda =\lambda \left( \delta \right) $ such that $\delta \exp
\left( \lambda \left( \delta \right) m\right) =\sqrt{\delta }.$ This means
that 
\begin{equation}
\lambda \left( \delta \right) =\ln \left( \frac{1}{\delta ^{1/\left(
2m\right) }}\right) .  \label{7.5}
\end{equation}%
The choice (\ref{7.5}) is possible since $\ln \left( \delta _{0}^{-1/\left(
2m\right) }\right) >\lambda _{1}$ and, therefore, $\lambda \left( \delta
\right) >\lambda _{1}$ for $\delta \in \left( 0,\delta _{0}\right) .$ Choose 
$\beta =e^{-\lambda \left( \delta \right) \varepsilon }.$ Hence, taking into
account (\ref{100}), we obtain 
\begin{equation}
\delta e^{\lambda \left( \delta \right) m}+e^{-\lambda \left( \delta \right)
\varepsilon /2}=\sqrt{\delta }+\delta ^{\varepsilon /\left( 4m\right) }\leq
2\delta ^{\theta }.  \label{7.6}
\end{equation}%
Thus, (\ref{7.4})-(\ref{7.6}) imply (\ref{2}). Next, (\ref{3}) is
established similarly with the part of the proof of Theorem 3.4 after (\ref%
{12.3}). $\square $

\section{Specifying equations}

\label{sec:8}

The scheme of section 3 is a general one and it can be applied to all three
main classes of Partial Differential Equations of the second order:
elliptic, parabolic and hyperbolic ones. Since the latter was explained in
detail in \cite{Kl1}, we only briefly specify these equations in this
section and formulate ill-posed Cauchy problems for them. So, Theorems 2.1,
3.1-3.5 can be reformulated for all three Cauchy problems considered in this
section.

\subsection{Quasilinear elliptic equation}

\label{sec:8.1}

We now rewrite the operator $A$ in (\ref{2.100}) as%
\begin{equation}
A_{ell}\left( u\right) =\dsum\limits_{i,j=1}^{n}a_{i,j}\left( x\right)
u_{x_{i}x_{j}}+A_{1}\left( x,\nabla u,u\right) ,x\in \Omega ,  \label{14.1}
\end{equation}%
\begin{equation}
A_{0,ell}\left( u\right) =\dsum\limits_{i,j=1}^{n}a_{i,j}\left( x\right)
u_{x_{i}x_{j}},  \label{14.2}
\end{equation}%
\begin{equation}
a_{i,j}\in C^{1}\left( \overline{\Omega }\right) ,  \label{14.3}
\end{equation}%
where $a_{i,j}\left( x\right) =a_{j,i}\left( x\right) ,\forall i,j=1,...,n$
and $A_{0}$ is the principal part of the operator $A_{ell}.$ Condition (\ref%
{300}) becomes now condition (\ref{14.3}). Also, we assume that condition (%
\ref{2.1002}) holds. The ellipticity of the operator $A_{0,ell}$ means that
there exist two constants $\mu _{1},\mu _{2}>0,\mu _{1}\leq \mu _{2}$ such
that%
\begin{equation}
\mu _{1}\left\vert \eta \right\vert ^{2}\leq
\dsum\limits_{i,j=1}^{n}a_{i,j}\left( x\right) \eta _{i}\eta _{j}\leq \mu
_{2}\left\vert \eta \right\vert ^{2},\forall x\in \overline{\Omega },\forall
\eta =\left( \eta _{1},...\eta _{n}\right) \in \mathbb{R}^{n}.  \label{14.4}
\end{equation}%
As above, let $\Gamma \subset \partial \Omega $ be the part of the boundary $%
\partial \Omega $, where the Cauchy data are given. Assume that the equation
of $\Gamma $ is 
\begin{equation*}
\Gamma =\left\{ x\in \mathbb{R}^{n}:x_{1}=p\left( \overline{x}\right) ,%
\overline{x}=\left( x_{2},...,x_{n}\right) \in \Gamma ^{\prime }\subset 
\mathbb{R}^{n-1}\right\}
\end{equation*}
and that the function $p\in C^{2}\left( \overline{\Gamma }^{\prime }\right)
. $ Here $\Gamma ^{\prime }\subset \mathbb{R}^{n-1}$ is a bounded domain.
Changing variables $x=\left( x_{1},\overline{x}\right) \Leftrightarrow
\left( x_{1}^{\prime },\overline{x}\right) ,$ where $x_{1}^{\prime
}=x_{1}-p\left( \overline{x}\right) $ and keeping the same notation for $%
x_{1}$ for brevity, we obtain that in new variables%
\begin{equation*}
\Gamma =\left\{ x\in \mathbb{R}^{n}:x_{1}=0,\overline{x}\in \Gamma ^{\prime
}\right\} .
\end{equation*}%
This change of variables does not affect the ellipticity property of the
operator $A_{0,ell}$. Let $X>0$ be a certain number. Thus, without any loss
of generality, we assume that 
\begin{equation}
\Omega \subset \left\{ x_{1}\in \left( 0,1/2\right) \right\} ,\text{ }\Gamma
=\left\{ x\in \mathbb{R}^{n}:x_{1}=0,\left\vert \overline{x}\right\vert
<X\right\} \subset \partial \Omega .  \label{14.5}
\end{equation}

\textbf{Cauchy Problem for the Quasilinear Elliptic Equation}. \emph{Suppose
that conditions (\ref{14.1})-(\ref{14.4}) hold. Find such a function }$u\in
H^{k}\left( \Omega \right) $\emph{\ that satisfies the equation }%
\begin{equation*}
A_{ell}\left( u\right) =0
\end{equation*}%
\emph{and has the following Cauchy data }$g_{0},g_{1}$ \emph{on }$\Gamma $%
\begin{equation*}
u\mid _{\Gamma }=g_{0}\left( \overline{x}\right) ,u_{x_{1}}\mid _{\Gamma
}=g_{1}\left( \overline{x}\right) .
\end{equation*}

Let $\rho \in \left( 0,1/2\right) $ be an arbitrary number. It is well known
that the CWF $\varphi _{\lambda }\left( x\right) $ for the operator $%
A_{0,ell}$ in this case can be chosen as 
\begin{equation}
\psi \left( x\right) =x_{1}+\frac{\left\vert \overline{x}\right\vert ^{2}}{%
X^{2}}+\rho ,\varphi _{\lambda }\left( x\right) =\exp \left[ \lambda \psi
^{-\nu }\left( x\right) \right] ,  \label{14.6}
\end{equation}%
see chapter 4 of \cite{LRS}. Here the number $\nu \geq \nu _{0},$ where $\nu
_{0}=\nu _{0}\left( \Omega ,n,\rho ,X,A_{0,ell}\right) >1$ is a certain
number depending only on listed parameters.

\subsection{Quasilinear parabolic equation}

\label{sec:8.2}

Since in this and next subsections we work with the space $\mathbb{R}%
^{n+1}=\left\{ \left( x,t\right) ,x\in \mathbb{R}^{n},t\in \mathbb{R}%
\right\} ,$ then we replace the above number $k$ with $k_{n+1}=\left[ \left(
n+1\right) /2\right] +2.$ Choose an arbitrary number $T=const.>0$ and denote 
$Q_{T}=\Omega \times \left( -T,T\right) .$ Let $L_{par}$ be the quasilinear
elliptic operator of the second order in $Q_{T},$ which we define the same
way as the operator $A_{ell}$ in (\ref{14.1})-(\ref{14.3}) with the only
difference that now its coefficients depend on both $x$ and $t,$ and also
the domain $\Omega $ is replaced with the domain $Q_{T}.$ Let $L_{0,par}$ be
the similarly defined principal part of the operator $L_{par},$ see (\ref%
{14.2}). Next, we define the quasilinear parabolic operator as $%
A_{par}=\partial _{t}-L_{par}$. The principal part of $A_{par}$ is $%
A_{0,par}=\partial _{t}-L_{0,par}.$ Thus, in $Q_{T}$ 
\begin{equation}
L_{par}\left( u\right) =\dsum\limits_{i,j=1}^{n}a_{i,j}\left( x,t\right)
u_{x_{i}x_{j}}+A_{1}\left( x,t,\nabla u,u\right) ,  \label{5.1}
\end{equation}%
\begin{equation}
A_{par}\left( u\right) =u_{t}-L_{par}\left( u\right) ,  \label{5.2}
\end{equation}%
\begin{equation}
A_{0,par}\left( u\right)
=u_{t}-L_{0,par}u=u_{t}-\dsum\limits_{i,j=1}^{n}a_{i,j}\left( x,t\right)
u_{x_{i}x_{j}},  \label{5.3}
\end{equation}%
\begin{equation}
a_{i,j}\in C^{1}\left( \overline{Q}_{T}\right) ,\text{ }  \label{5.4}
\end{equation}%
\begin{equation}
A_{1}\left( x,t,y\right) \in C^{3}\left\{ \left( x,t,y\right) :\left(
x,t\right) \in \overline{Q}_{T},y\in \mathbb{R}^{n+1}\right\} ,  \label{5.5}
\end{equation}%
\begin{equation}
\mu _{1}\left\vert \eta \right\vert ^{2}\leq
\dsum\limits_{i,j=1}^{n}a_{i,j}\left( x,t\right) \eta _{i}\eta _{j}\leq \mu
_{2}\left\vert \eta \right\vert ^{2},\forall \left( x,t\right) \in \overline{%
Q}_{T},\forall \eta =\left( \eta _{1},...\eta _{n}\right) \in \mathbb{R}^{n}.
\label{5.6}
\end{equation}%
Let the domain $\Omega $ and the hypersurface $\Gamma \subset \partial
\Omega $ be the same as in (\ref{14.5}). Denote $\Gamma _{T}=\Gamma \times
\left( -T,T\right) .$ Consider the quasilinear parabolic equation 
\begin{equation}
A_{par}\left( u\right) =u_{t}-L_{par}\left( u\right) =0\text{ \ in }Q_{T}.%
\text{ }  \label{5.7}
\end{equation}

\textbf{Cauchy Problem with the Lateral Data for the Quasilinear Parabolic
Equation.} \emph{Assume that conditions (\ref{5.1})-(\ref{5.6}) hold. Find
such a function }$u\in H^{k_{n+1}}\left( Q_{T}\right) $\emph{\ which
satisfies equation (\ref{5.7}) and has the following lateral Cauchy data }$%
g_{0},g_{1}$ \emph{at }$\Gamma _{T}$%
\begin{equation}
u\mid _{\Gamma _{T}}=g_{0}\left( \overline{x},t\right) ,u_{x_{1}}\mid
_{\Gamma _{T}}=g_{1}\left( \overline{x},t\right) .  \label{5.70}
\end{equation}

The CWF $\varphi _{\lambda }\left( x,t\right) $ for the operator $A_{0,par}$
is introduced similarly with (\ref{14.6}), see chapter 4 of \cite{LRS} 
\begin{equation}
\psi \left( x,t\right) =x_{1}+\frac{\left\vert \overline{x}\right\vert ^{2}}{%
X^{2}}+\frac{t^{2}}{T^{2}}+\rho ,\text{ }\varphi _{\lambda }\left(
x,t\right) =\exp \left( \lambda \psi ^{-\nu }\right) .  \label{14.7}
\end{equation}%
Here the number $\nu \geq \nu _{0},$ where $\nu _{0}=\nu _{0}\left( \Omega
,n,\rho ,X,T,L_{0,par}\right) >1$ is a certain number depending only on
listed parameters.

This CWF works perfectly for the case when the lateral Cauchy data are given
on a part $\Gamma \neq \partial \Omega $ of the boundary of the domain $%
\Omega $. Suppose now that $\Omega $ is a ball, $\Omega =\left\{ \left\vert
x\right\vert <B\right\} $ for a constant $B>0.$ Suppose that the data are
given on the entire boundary $\partial \Omega =\left\{ \left\vert
x\right\vert =B\right\} .$ Then one can use \cite{Y} 
\begin{equation}
\varphi _{\lambda }\left( x,t\right) =\exp \left[ \lambda \left( \left\vert
x\right\vert ^{2}-t^{2}\right) \right] .  \label{14.8}
\end{equation}

\subsection{Quasilinear hyperbolic equation}

\label{sec:8.3}

In this subsection, notations for the time cylinder $Q_{T}$ are the same as
ones in subsection 8.2. We assume here that $\partial \Omega \in C^{\infty
}. $ Denote $S_{T}=\partial \Omega \times \left( -T,T\right) .$ Consider two
numbers $a_{l},a_{u}>0$ such that $a_{l}<a_{u}.$ For $x\in \Omega ,$ let the
function $a\left( x\right) $ satisfy the following conditions 
\begin{equation}
a\left( x\right) \in \left[ a_{l},a_{u}\right] ,a\in C^{1}\left( \overline{%
\Omega }\right) .  \label{6.1}
\end{equation}%
In addition, we assume that there exists a point $x_{0}\in \mathbb{R}%
^{n}\diagdown \overline{\Omega }$ such that%
\begin{equation}
\left( \nabla a\left( x\right) ,x-x_{0}\right) \geq 0,\text{ }\forall x\in 
\overline{\Omega }.  \label{6.2}
\end{equation}%
In particular, if $a\left( x\right) \equiv 1,$ then (\ref{6.2}) holds for
any $x_{0}\in \mathbb{R}^{n}.$ We need inequality (\ref{6.2}) for the
validity of the Carleman estimate. Assume that the function $A_{1}$
satisfies condition (\ref{5.5}). Consider the quasilinear hyperbolic
equation in the time cylinder $\Omega _{T}$ with the lateral Cauchy data $%
g_{0}\left( x,t\right) ,g_{1}\left( x,t\right) $ at $S_{T},$ 
\begin{eqnarray}
L_{hyp}\left( u\right) &=&a\left( x\right) u_{tt}-\Delta u-A_{1}\left(
x,t,\nabla u,u\right) =0\text{ in }Q_{T},  \label{6.3} \\
u &\mid &_{S_{T}}=g_{0}\left( x,t\right) ,\partial _{n}u\mid
_{S_{T}}=g_{1}\left( x,t\right) .  \label{6.4}
\end{eqnarray}

\textbf{Cauchy Problem with the Lateral Data for the Hyperbolic Equation (%
\ref{6.3})}. \emph{Find such a function }$u\in H^{k_{n+1}}\left(
Q_{T}\right) $\emph{\ which satisfies conditions (\ref{6.3}) and (\ref{6.4}).%
}

Let the number $\eta \in \left( 0,1\right) .$ Define functions $\psi \left(
x,t\right) $ and $\varphi _{\lambda }\left( x,t\right) $ as%
\begin{equation}
\psi \left( x,t\right) =\left\vert x-x_{0}\right\vert ^{2}-\eta t^{2},\text{ 
}\varphi _{\lambda }\left( x,t\right) =\exp \left( \lambda \psi \left(
x,t\right) \right) .  \label{6.5}
\end{equation}%
Denote $L_{0,hyp}\left( u\right) =a\left( x\right) u_{tt}-\Delta u.$ Given
conditions (\ref{6.1}), (\ref{6.2}), the Carleman estimate for the operator $%
L_{0,hyp}$ holds with the CWF $\varphi _{\lambda }$ from (\ref{6.5}), see
Theorem 1.10.2 in the book \cite{BK}.

\section{Numerical Study}

\label{sec:9}

In this section, we study numerically a 1-D analog of the ill-posed Cauchy
problem (\ref{5.7}), (\ref{5.70}) for the parabolic equation. The numerical
study of this section is similar with the one of \cite{KKLY}. There are
important differences, however. First, following the concept of Cauchy
Problem 2, we obtain zero Dirichlet and Neumann boundary conditions on one
edge of the interval, where the lateral Cauchy data are given. Second, the
specific formulas for the quasilinear part $S\left( u\right) $ of the
parabolic operator considered below are different from ones of \cite{KKLY}.
Still, because of the above analogy, our description below is rather brief.
We refer to \cite{KKLY} for more details.

\subsection{The forward problem}

\label{sec:9.1}

Here $T=1/2$ and 
\begin{equation*}
Q_{1/2}=\left\{ \left( x,t\right) :x\in \left( 0,1\right) ,t\in \left(
-1/2,1/2\right) \right\} .
\end{equation*}%
We consider the following forward problem:%
\begin{equation}
u_{t}=u_{xx}+S\left( u\right) +G\left( x,t\right) ,\left( x,t\right) \in
Q_{T},  \label{15.1}
\end{equation}%
\begin{equation}
u\left( x,-1/2\right) =f\left( x\right) ,  \label{15.2}
\end{equation}%
\begin{equation}
u\left( 0,t\right) =g\left( t\right) ,u\left( 1,t\right) =p\left( t\right) .
\label{15.3}
\end{equation}%
Our specific functions in (\ref{15.1})-(\ref{15.3}) are:%
\begin{equation}
S_{1}\left( u\right) =10\cos \left( u+x+2t\right) ,S_{2}\left( u\right) =10%
\frac{u^{2}}{1+u^{2}},  \label{15.4}
\end{equation}%
\begin{equation}
G\left( x,t\right) =10\sin \left[ 100\left( \left( x-0.5\right)
^{2}+t^{2}\right) \right] ,  \label{15.5}
\end{equation}%
\begin{equation}
f\left( x\right) =10\left( x-x^{2}\right) ,  \label{15.6}
\end{equation}%
\begin{equation}
g\left( t\right) =10\sin \left[ 10\left( t-0.5\right) \left( t+0.5\right) %
\right] ,p\left( t\right) =\sin \left[ 10\left( t+0.5\right) \right] .
\label{15.7}
\end{equation}%
Thus, due to the presence of the multiplier $10$ in (\ref{15.4}), the
influence of the nonlinear term $S\left( u\right) $ on the solution $u$ of
the problem (\ref{15.1})-(\ref{15.3}) is significant.

We use the FDM to solve the forward problem (\ref{15.1})-(\ref{15.3})
numerically. Introduce the uniform mesh in the domain $\Omega _{T},$%
\begin{equation*}
\overline{M}=\left\{ \left( x_{i},t_{j}\right) :x_{i}=ih,t_{j}=-\frac{1}{2}%
+j\tau ,i\in \left[ 0,N\right) ,j\in \left[ 0,M\right) \right\} ,
\end{equation*}%
where $h=1/N$ and $\tau =1/M$ are grid step sizes in $x$ and $t$ directions
respectively. For generic functions $f^{\left( 1\right) }\left( x,t\right)
,f^{\left( 2\right) }\left( x\right) ,f^{\left( 3\right) }\left( t\right) $
denote $f_{ij}^{\left( 1\right) }=f^{\left( 1\right) }\left(
x_{i},t_{j}\right) ,f_{i}^{\left( 2\right) }=f^{\left( 2\right) }\left(
x_{i}\right) ,f_{j}^{\left( 3\right) }=f^{\left( 3\right) }\left(
t_{j}\right) .$ Let $\varphi _{ij}=S\left( u_{ij}\right) +G_{ij}.$ We have
solved the forward problem (\ref{15.1})-(\ref{15.3}) using the implicit
finite difference scheme,%
\begin{equation*}
\frac{u_{ij+1}-u_{ij}}{\tau }=\frac{1}{h^{2}}\left(
u_{i-1j+1}-2u_{ij+1}+u_{i+1j+1}\right) +\varphi _{ij},i\in \left[
1,N-1\right) ,j\in \left[ 0,M-1\right) ,
\end{equation*}%
\begin{equation*}
u_{i0}=f_{i},\text{ }u_{0j}=g_{j},u_{Nj}=p_{j},i\in \left[ 0,N\right) ,j\in %
\left[ 0,M\right) .
\end{equation*}%
In all our numerical tests we have used $M=32,N=128$. Even though these
numbers are the same both for the solution of the forward and inverse
problems, the \textquotedblleft inverse crime" was not committed since we
have used noisy data and since we have used the minimization of a functional
rather than solving a forward problem again.\emph{\ }

Thus, solving the forward problem (\ref{15.1})-(\ref{15.3}) with the input
functions (\ref{15.4})-(\ref{15.7}), we have computed the function $%
q_{comp}\left( t\right) ,$%
\begin{equation}
u_{x}\left( 1,t\right) =q_{comp}\left( t\right) ,t\in \left( -\frac{1}{2},%
\frac{1}{2}\right) .  \label{15.8}
\end{equation}

\subsection{The ill-posed Cauchy problem and noisy data}

\label{sec:9.2}

Our interest in this section is to solve numerically the following \ Cauchy
problem:

\textbf{1-D Cauchy Problem}. \emph{Suppose that in (\ref{15.1})-(\ref{15.3})
functions }$f\left( x\right) $\emph{\ and }$g\left( t\right) $\emph{\ are
unknown whereas the functions }$G\left( x,t\right) ,$\emph{\ }$p\left(
t\right) $\emph{\ and }$S\left( u\right) $\emph{\ are known. Let in the data
simulation process functions }$F,S,f,g,p$\emph{\ are the same as in (\ref%
{15.4})-(\ref{15.7}). Determine the function }$u\left( x,t\right) $\emph{\
in at least a subdomain of the time cylinder }$Q_{1/2},$\emph{\ assuming
that the function }$q_{comp}\left( t\right) $\emph{\ in (\ref{15.8}) is
known.}

We have introduced 5\% level of random noise in the data. Let $\sigma \in %
\left[ -1,1\right] $ be the random variable representing the white noise.
Let $p^{\left( m\right) }=\max_{j}\left\vert p_{j}\right\vert $ and $%
q^{\left( m\right) }=\max_{j}\left\vert q_{comp,j}\right\vert .$ Then the
noisy data, which we have used, were%
\begin{equation}
\widetilde{u}_{N-1j}=p_{j}+0.05p^{\left( m\right) }\sigma _{j},\widetilde{u}%
_{N-2j}=p_{j}-h\left( q_{comp,j}+0.05q^{\left( m\right) }\sigma _{j}\right) .
\label{5.16}
\end{equation}%
Below we use functions $p^{\prime }\left( t\right) $ and $q^{\prime }\left(
t\right) .$ We have calculated derivatives $p^{\prime }\left( t\right)
,q^{\prime }\left( t\right) $ of noisy functions via finite differences.
Even though the differentiation of noisy functions is an ill-posed problem,
we have not observed instabilities in our case. A more detailed study of
this topic is outside of the scope of the current publication.

\subsection{Specifying the functional $J_{\protect\lambda ,\protect\beta }$}

\label{sec:9.3}

We introduce the function $F(x,t)$ as%
\begin{equation*}
F(x,t)=p\left( t\right) +\left( x-1\right) q\left( t\right) .
\end{equation*}%
Let $v=u-F.$ Then $v\left( 1,t\right) =v_{x}\left( 1,t\right) =0$ and 
\begin{equation*}
A\left( v+F\right) =v_{t}-v_{xx}-S\left( v+F\right) -G\left( x,t\right) + 
\left[ p^{\prime }\left( t\right) +\left( x-1\right) q^{\prime }\left(
t\right) \right] .
\end{equation*}%
By (\ref{2.1005})\emph{\ }the functional $J_{\lambda ,\beta }\left(
v,F\right) $ becomes%
\begin{equation}
J_{\lambda ,\beta }\left( v,F\right)
=\dint\limits_{-1/2}^{1/2}\dint\limits_{0}^{1}\left[ A\left( v+F\right) %
\right] ^{2}\varphi _{\lambda }^{2}dxdt+\beta \left\Vert v\right\Vert
_{H^{2}\left( Q_{1/2}\right) }^{2}.  \label{15.9}
\end{equation}%
Here we use the CWF $\varphi _{\lambda }\left( x,t\right) $ which is given
in (\ref{14.8}). The reason of this is that the rate of change of the CWF of
(\ref{14.7}) is too large due to the presence of two large parameters $%
\lambda $ and $\nu $ in (\ref{14.7}) In (\ref{15.9}) we do not use the
multiplier $e^{-2\lambda \left( \alpha +\varepsilon \right) },$ which was
present in the original version (\ref{2.1005}). Indeed, we have used that
multiplier in order to allow the parameter $\beta $ to be less than $1$.
However, we have observed in our computations that the accuracy of results
does not change much for $\beta $ varying in a large interval. In all our
numerical experiments below $\beta =0.00063.$ The norm $\left\Vert
v\right\Vert _{H^{2}\left( Q_{T}\right) }$ is taken instead of $\left\Vert
v\right\Vert _{H^{3}\left( Q_{1/2}\right) }$ due to the convenience of
computations. Note that since we do not use too many grid points when
discretizing the functional $J_{\lambda ,\beta }\left( v,F\right) ,$ then
these two norms are basically equivalent in our computations, since all
norms are equivalent in a finite dimensional space.

\subsection{Minimization of $J_{\protect\lambda ,\protect\beta }\left(
v,F\right) $}

\label{sec:9.4}

To minimize the functional (\ref{15.9}), we have attempted first to use the
gradient projection method, as it was done in the above theoretical part.
However, we have observed in our computations that just the conjugate
gradient method (GCM) with the starting point $v_{0}\equiv 0$ works well and
much more rapidly. So, our results below are obtained via the GCM. We have
written the functional $J_{\lambda ,\beta }\left( v,F\right) $ in the
discrete form $\overline{J}_{\lambda ,\beta }\left( v,F\right) $ using
finite differences. Next, we have minimized the functional $\overline{J}%
_{\lambda ,\beta }\left( v,F\right) $ with respect to the values $v_{ij}$ of
the discrete function $v$ at the grid points. Hence, we have calculated
derivatives $\partial _{v_{ij}}\overline{J}_{\lambda ,\beta }\left(
v,F\right) $ via explicit formulas. The method of the calculation of these
derivatives is described in \cite{KKLY}.

Normally, for a quadratic functional the GCM reaches the minimum of this
functional after $M\cdot N$ gradient steps with the automatic step choice.
However, our computational experience tells us that we can obtain a better
accuracy if using a small constant step in the GCM and a large number of
iterations. Thus, we have used the step size $\gamma =10^{-8}$ and 10,000
iterations of the GCM. It took 0.5 minutes of CPU Intel Core i7 to do these
iterations.

\subsection{Results}

\label{sec:9.5}

Let $v\left( x,t\right) $ be the numerical solution of the forward problem (%
\ref{15.1})-(\ref{15.3}). Let $v_{\lambda \beta }\left( x,t\right) $ be the
minimizer of the functional $\overline{J}_{\lambda ,\beta }\left( v,F\right) 
$ which we have found via the GCM. Of course, $v\left( x,t\right) $ and $%
v_{\lambda \beta }\left( x,t\right) $ here are discrete functions defined on
the above grid and norms used below are discrete norms. Recall that $u\left(
x,t\right) =v\left( x,t\right) +F\left( x,t\right) .$ Hence, denote $%
u_{\lambda \beta }\left( x,t\right) =v_{\lambda \beta }\left( x,t\right)
+F\left( x,t\right) .$ For each $x$ of our grid we define the
\textquotedblleft line error" $E\left( x\right) $ as%
\begin{equation}
E\left( x\right) =\frac{\left\Vert u_{\lambda \beta }\left( x,t\right)
-u\left( x,t\right) \right\Vert _{L_{2}\left( -1/2,1/2\right) }}{\left\Vert
u\left( x,t\right) \right\Vert _{L_{2}\left( -1/2,1/2\right) }}.
\label{15.10}
\end{equation}%
We evaluate how the line error changes with the change of $x$, i.e. how the
computational error changes when the point $x$ moves away from the edge $x=1$
where the lateral Cauchy data are given. Naturally, it is anticipated that
the function $E\left( x\right) $ should be decreasing.

\begin{figure}[ht!]
\begin{center}
\begin{tabular}{c}
\includegraphics[width=0.9\linewidth,clip=]{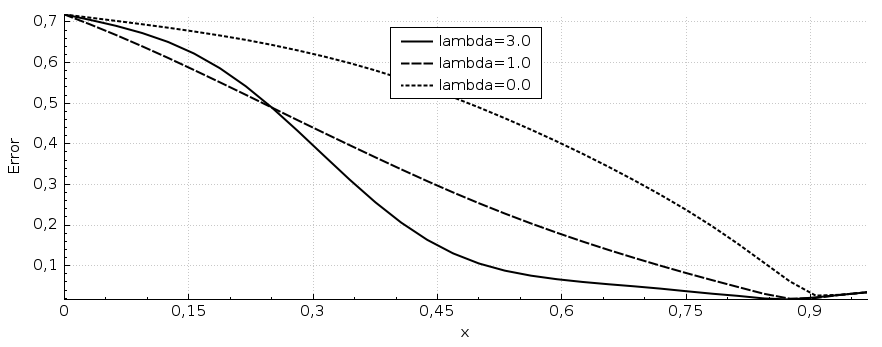} \\ 
a) \\ 
\includegraphics[width=0.9\linewidth,clip=]{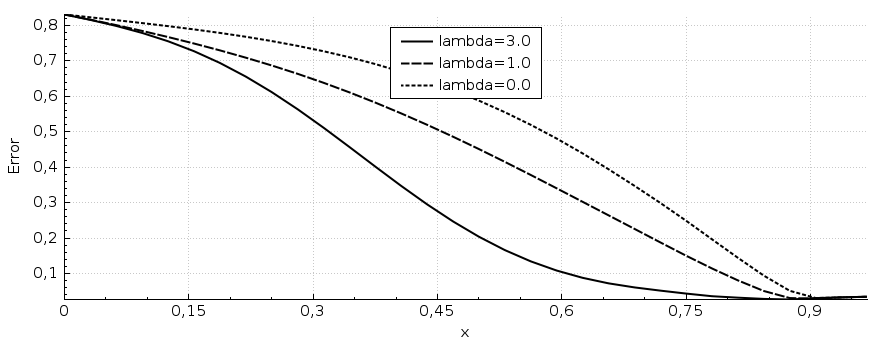} \\ 
b)%
\end{tabular}%
\end{center}
\caption{ \emph{Distribution of error along the x-axis. a) }$\protect\lambda %
=0,1,3$\emph{\ and }$S\left( u\right) =\cos\left( u + x +2t \right) .$\emph{%
\ b) }$\protect\lambda =0,1,3$\emph{\ \ and }$S\left( u\right) = u^2/(1+u^2)
.$\emph{\ Thus, the presence of the CWF in the functional (\protect\ref{15.9}%
) significantly improves the accuracy of the solution. One can observe that
a rather accurate reconstruction is obtained on the interval $[0.6,1]$,
also, see Remark 9.1.}}
\end{figure}

\textbf{Remark 9.1}. It is clear, intuitively at least, that the further a
point $x\in \left( 0,1\right) $ is from the point $x=1$ where the lateral
Cauchy data are given, the less accuracy of solution at this point one
should anticipate. So, we observe in graphs of line errors on Figures
1a)-4a) that the accuracy of the calculated solutions for $x\in \left(
0,0.6\right) $ is not as good as this accuracy for $x\in \left[ 0.6,1\right]
.$ This is why we graph below only line errors and functions $u_{\lambda
\beta }\left( 0.6,t\right) ,$ superimposed with $u\left( 0.6,t\right) .$

In the case of Figures 1 and 2 the starting function for iterations of the
GCM was $v_{0}\equiv 0.$ We have tested three values of the parameter $%
\lambda :\lambda =0,1,3$ in (\ref{14.8}). We have found that $\lambda =3$ is
the best choice for those problems which we have studied. This is also clear
from Figures 1. Note that the case $\lambda =0$ provides a poor accuracy.

As one can see on Figures 1, the line error at $x=0.6$ is between about 6\%
\ and 10\% for $\lambda =3$. Thus, we superimpose graphs of functions $%
u_{\lambda \beta }\left( 0.6,t\right) $ with graphs of functions $u\left(
0.6,t\right) $ (see Remark 9.1)$.$ Corresponding graphs are displayed on
Figures 2. One can observe again that the computational accuracy with $%
\lambda =3$ is the best and that the accuracy with $\lambda =0$ is poor.
Thus, we observe again that the presence of the CWF in the functional (\ref%
{15.9}) significantly improves the accuracy of the solution. On the other
hand, the accuracy at $t\approx \pm 1/2$ is not good on Figures 2. We
explain this by the fact that Theorem 3.5 guarantees a good accuracy only in
a subdomain $\Omega _{\alpha +2\varepsilon }$ of the domain $\Omega $ rather
than in the entire domain $\Omega .$ The latter can be reformulated for our
specific domain $Q_{1/2}$ \cite{KKLY}.

\begin{figure}[hb!]
\begin{center}
\begin{tabular}{c}
\includegraphics[width=0.9\linewidth,clip=]{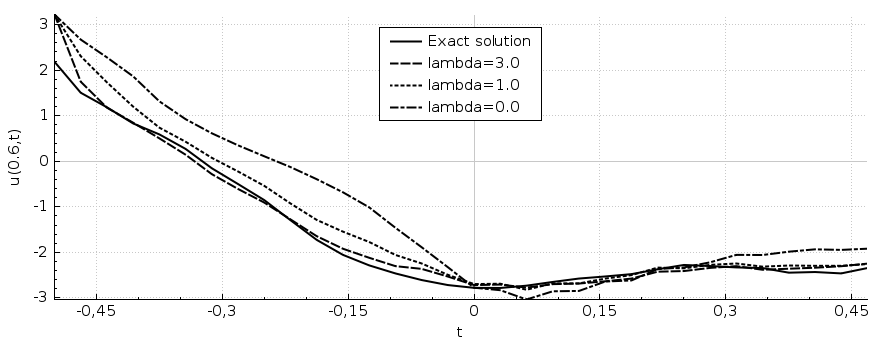} \\ 
a) \\ 
\includegraphics[width=0.9\linewidth,clip=]{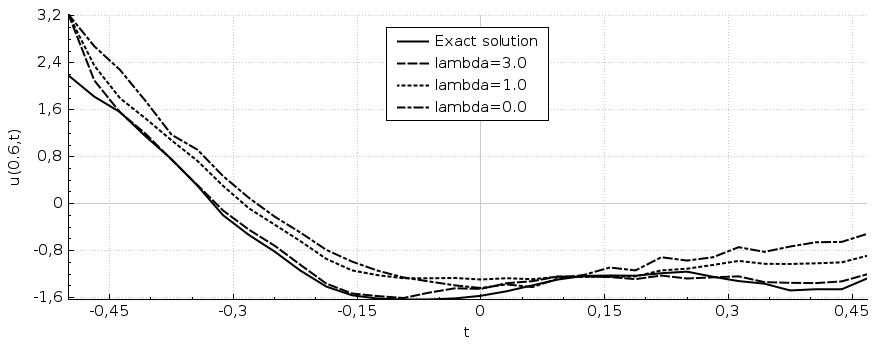} \\ 
b)%
\end{tabular}%
\end{center}
\caption{ \emph{Superimposed graphs of functions }$u_{0\protect\beta }\left(
0.6,t\right) ,u_{1\protect\beta }\left( 0.6,t\right),u_{3\protect\beta %
}\left( 0.6,t\right) $\emph{\ and }$u\left( 0.6,t\right) .$\emph{\ a) }$%
S\left( u\right) =\cos\left( u + x +2t \right) .$\emph{\ b) }$S\left(
u\right) =u^2/(1+u^2) .$\emph{\ Observe that the presence of the CWF with }$%
\protect\lambda =3$\emph{\ significantly improves the accuracy of the
solution.} }
\label{fig:slices}
\end{figure}

To see how the starting function of the GCM affects the accuracy of our
results, we took $S\left( u\right) =S_{1}\left( u\right) $ and have tested
three starting functions for the GCM: $v_{0}\left( x,t\right) =0,v_{0}\left(
x,t\right) =\left( x-1\right) ^{2}\left( t+1\right) $ and $v_{0}\left(
x,t\right) =\left( \sin \left( x-1\right) \right) ^{2}t^{2}$. Hence, for any
of these three functions $v_{0}\left( x,t\right) $ we have $v_{0}\left(
1,t\right) =\partial _{x}v_{0}\left( 1,t\right) =0.$ Graphs of Figure 3a)
displays superimposed line errors and Figure 3b) displays functions $%
u_{3\beta }\left( 0.6,t\right) $ and $u\left( 0.6,t\right) $ for these three
cases (see Remark 9.1). One can see that for $x\in \left[ 0.6,1\right] $
results depend only very insignificantly on the starting point of the GCM:
just as it was predicted by Theorems 3.3 and 3.5, also see Remark 9.1.

\begin{figure}[hb!]
\begin{center}
\begin{tabular}{c}
\includegraphics[width=0.9\linewidth,clip=]{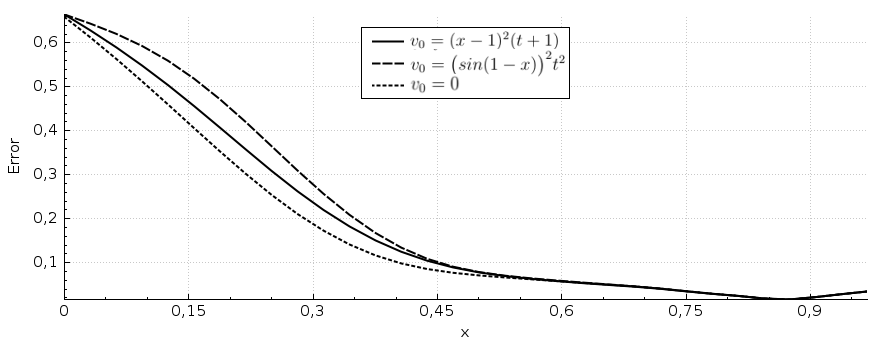} \\ 
a) \\ 
\includegraphics[width=0.9\linewidth,clip=]{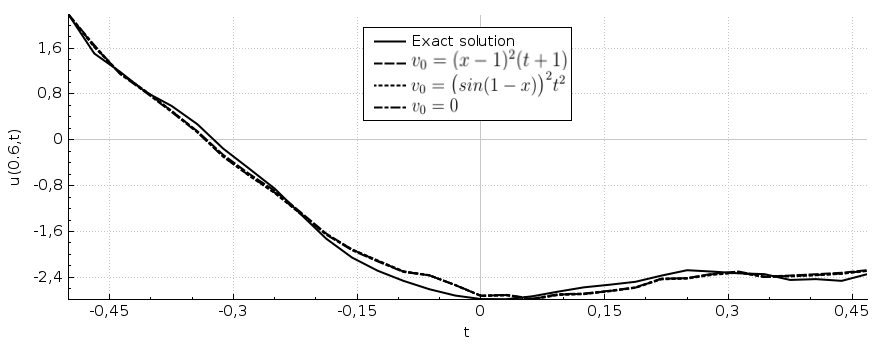} \\ 
b)%
\end{tabular}%
\end{center}
\caption{ \emph{The influence of the choice of the starting function }$v_{0}$
\emph{of the GCM. We have tested three starting functions: }$%
v_{0}=0,v_{0}=(x-1)^2(t+1)$\emph{\ and }$v_{0}=\big(sin(1-x)\big)^2t^2$\emph{%
\ We took }$\protect\lambda =3,S\left( u\right) =10\cos \left( u+x+2t\right)
.$\emph{\ a) Superimposed line errors. b) Superimposed functions }$u\left(
0.6,t\right) $\emph{\ and }$u_{3\protect\beta }\left( 0.6,t\right) .$
\emph{These tests demonstrate that for }$x\in \left[ 0.6,1\right] $\emph{\
our solution depends only very insignificantly from the choice of the
starting function }$v_{0}$\emph{\ of the GCM: just as it was predicted by
Theorems 3.3 and 3.5. Also, see Remark 9.1.} }
\label{fig:initial_cond_influence}
\end{figure}

It is interesting to see how the presence of the CWF affects the linear
case. In this case the above method with $\lambda =0$ becomes the
Quasi-Reversibility Method \cite{LL}. In the recent survey of the second
author convergence rates of regularized solutions were estimated for this
method \cite{KQ}. The existence of regularized solutions, i.e. minimizers,
was also established in \cite{KQ}. So, we have tested the case when 
\begin{equation*}
S\left( u\right) \equiv 0
\end{equation*}%
in (\ref{15.1}), while functions $G\left( x,t\right) ,f\left( x\right)
,g\left( t\right) $ and $p\left( t\right) $ are the same as in (\ref{15.5})-(%
\ref{15.7}). Results for $\lambda =0,1,3$ are presented on Figures 4. One
can observe that even in the linear case the presence of the CWF
significantly improves the computational accuracy for $x\in \left[ 0.6,1%
\right] $.

\begin{figure}[ht!]
\begin{center}
\begin{tabular}{c}
\includegraphics[width=0.9\linewidth,clip=]{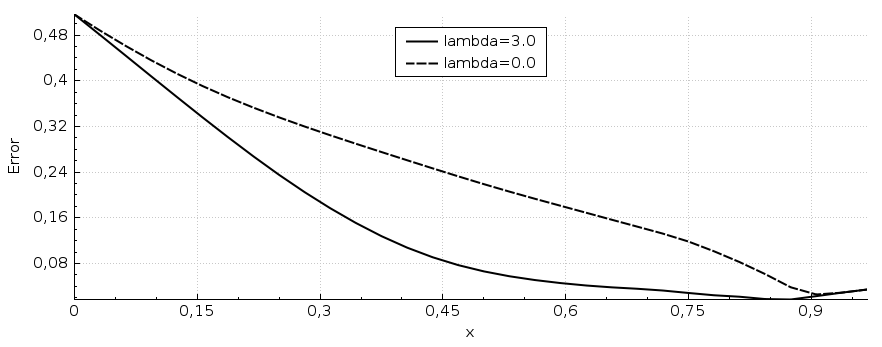} \\ 
a) \\ 
\includegraphics[width=0.9\linewidth,clip=]{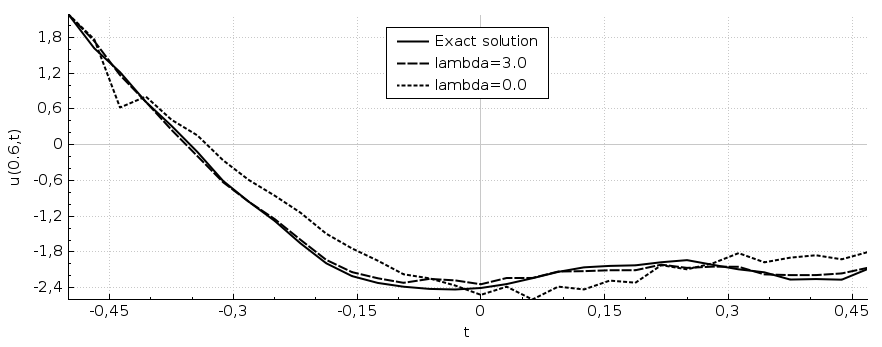} \\ 
b)%
\end{tabular}%
\end{center}
\caption{ \emph{The linear case: when }$S\left( u\right) \equiv 0$ \emph{in (%
\protect\ref{15.1}). a) Line errors. b) Functions }$u\left( 0.6,t\right) $%
\emph{\ and }$u_{3\protect\beta }\left( 0.6,t\right) .$
\emph{One can observe that the presence of the CWF significantly improves
the accuracy of the solution for }$x\in \left[ 0.6,1\right] $\emph{\ even in
the linear case. Also, see Remark 9.1.} }

\label{fig:linear_case}
\end{figure}

\section{Summary}

\label{sec:10}

In this paper a new element is introduced in the theory, which was
previously developed in \cite{BKconv,Klib97,Kpar,KNT,Kl1,KK,KKLY}. More
precisely, we presented some facts of the convex analysis in section 2.
Next, using, as an example, a general ill-posed Cauchy problem for a
quasilinear PDE of the second order, we have shown that these facts ensure
the existence of the minimizer of a weighted Tikhonov-like functional on any
closed ball in a reasonable Hilbert space. This functional is strictly
convex on that ball. The strict convexity is due to the presence of the
Carleman Weight Function. Next, we have specified PDEs of the second order
for which this construction works. We have also pointed out in section 1
that similar results are valid for coefficient inverse problems, which were
considered in \cite{BKconv,Klib97,Kpar,KNT,KK}. On the other hand, the
existence of the minimizer was assumed rather than proved in those previous
publication.

It is because of the material of section 2, that we have proved Theorem 3.1
and have significantly changed the methods of proofs of Theorems 3.2 and
3.5, as compared with \cite{Kl1}. Next, we have specified quasilinear PDEs
of the second order for which our technique works.

In addition, we have presented some numerical results for the side Cauchy
problem for a 1-D parabolic PDEs. These results indicate that the presence
of the CWF significantly improves the accuracy of the solution. Furthermore,
this is also true even in the linear case. It was also demonstrated
numerically that for $x\in \left[ 0.6,1\right] $ our resulting solution
depends on the starting function for the GCM only very insignificantly: just
as it was predicted by our theory, also see Remark 9.1.

\begin{center}
\textbf{Acknowledgments}
\end{center}

The work of A.B. Bakushinskii was supported by grants 15-01-00026 and
16-01-00039 of the Russian Foundation for Basic Research. The work of M.V.
Klibanov was supported by the US Army Research Laboratory and the US Army
Research Office grant W911NF-15-1-0233 as well as by the Office of Naval
Research grant N00014-15-1-2330.


\begin{thebibliography}{99}
\bibitem{Koz} S. Avdonin, V. Kozlov, D. Maxwell and M. Truffer, Iterative
methods for solving a nonlinear boundary inverse problem in glaciology, 
\emph{J. Inverse and Ill-Posed Problems}, 17, 239-258, 2009.

\bibitem{Alif1} O.M. Alifanov, \emph{Inverse Heat Conduction Problems},
Springer, New York, 1994.

\bibitem{Alif2} O.M. Alifanov, E.A.\ Artukhin and S.V.\ Rumyantcev, \emph{%
Extreme Methods for Solving Ill-Posed Problems with Applications to Inverse
Heat Transfer Problems, Begell House}, New York, 1995.

\bibitem{BK} L. Beilina and M.V. Klibanov, \emph{Approximate Global
Convergence and Adaptivity for Coefficient Inverse Problems}, Springer, New
York, 2012.

\bibitem{BKconv} L. Beilina and M.V. Klibanov, Globally strongly convex cost
functional for a coefficient inverse problem, \emph{Nonlinear Analysis:\
Real World Applications}, 22, 272-278, 2015.

\bibitem{BukhK} A.L. Bukhgeim and M.V. Klibanov, Uniqueness in the large of
a class of multidimensional inverse problems, \emph{Soviet Mathematics
Doklady}, 17, 244-247, 1981.

\bibitem{Bukh} A.L. Bukhgeim, \emph{Introduction to the Theory of Inverse
Problems}, VSP, Utrecht, 2000.

\bibitem{Col} J. Colinge and J. Rappaz, A strongly nonlinear problem arising
in glaceology, \emph{Mathematical Modelling and Numerical Analysis, M2AN},
33, 395-406, 1999.

\bibitem{Is} V. Isakov, \emph{Inverse Problems for Partial Differential
Equations}, 2nd Edition, Springer, New York, 2006.

\bibitem{K92} M.~V.~Klibanov, Inverse problems and Carleman estimates, \emph{%
Inverse Problems}, 8, 575--596, 1992.\emph{ems}, 7, 577-596, 1991.

\bibitem{Klib97} M.V. Klibanov, Global convexity in a three-dimensional
inverse acoustic problem, \emph{SIAM J.\ Mathematical Analysis}, 28,
1371-1388, 1997.

\bibitem{Kpar} M.V. Klibanov, Global convexity in diffusion tomography, 
\emph{Nonlinear World}, 4, 247-265, 1997.

\bibitem{KT} M.V. Klibanov and A. Timonov, \emph{Carleman Estimates for
Coefficient Inverse Problems and Numerical Applications}, VSP, Utrecht, 2004.

\bibitem{Ksurvey} M.V. Klibanov, Carleman estimates for global uniqueness,
stability and numerical methods for coefficient inverse problems, \emph{J.
Inverse and Ill-Posed Problems}, 21, 477-560, 2013.

\bibitem{KQ} M.V. Klibanov, Carleman estimates for the regularization of
ill-posed Cauchy problems, \emph{Applied Numerical Mathematics}, 94, 46-740,
2015.

\bibitem{KNT} M.V. Klibanov and N.T. Th\`{a}nh, Recovering of dielectric
constants of explosives via a globally strictly convex cost functional, 
\emph{SIAM J. Applied Mathematics}, 75, 518-537, 2015.

\bibitem{Kl1} M.V. Klibanov, Carleman weight functions for solving ill-posed
Cauchy problems for quasilinear PDEs, \emph{Inverse Problems}, 31, 125007,
2015.

\bibitem{KK} M.V. Klibanov and V.G. Kamburg, Globally strictly convex cost
functional for an inverse parabolic problem, \emph{Mathematical Methods in
the Applied Sciences, }39, 930-940, 2016.

\bibitem{KKLY} M.V. Klibanov, N.A. Koshev, J. Li and A.G. Yagola, Numerical
solution of an ill-posed Cauchy problem for a quasilinear parabolic equation
using a Carleman weight function, \emph{arxiv}: 1603.00848, 2016; accepted for 
publication in J. Inverse and Ill-Posed Probles.

\bibitem{LL} R. Lattes, J.-L. Lions, The Method of Quasireversibility:
Applications to Partial Differential Equations, Elsevier, New York, 1969.

\bibitem{LRS} M.M. Lavrentiev, V.G. Romanov and S.P. Shishatskii, \emph{%
Ill-Posed Problems of Mathematical Physics and Analysis}, AMS,\ Providence,
RI,\ 1986.

\bibitem{Scales} J. A. Scales, M.L. Smith and T.L. Fisher, Global
optimization methods for multimodal inverse problems \emph{J. Computational
Physics}, 103, 258-268, 1992.

\bibitem{T} A.N. Tikhonov, A.V. Goncharsky, V.V. Stepanov and A.G. Yagola, 
\emph{Numerical Methods for the Solution of Ill-Posed Problems}, Kluwer,
London, 1995.

\bibitem{Tit} V.N.Titarenko and A.G.Yagola. Cauchy problems for Laplace
equation on compact sets, \emph{Inverse Problems in Engineering}, 10,
235-254, 2002.

\bibitem{Trig} R. Triggiani and P.F. Yao, Carleman estimates with no lower
order terms for general Riemannian wave equations. Global uniqueness and
observability in one shot, \emph{Applied Mathematics and Optimization}, 46,
331-375, 2002.

\bibitem{Vas} F.P. Vasiliev, \emph{Numerical Methods of Solutions of
Extremal Problems}, Moscow, Nauka, 1989.

\bibitem{Y} M. Yamamoto, Carleman estimates for parabolic equations, \emph{%
Inverse Problems}, 25, 123013, 2009.
\end{thebibliography}
\end{document}